\documentclass[12pt]{amsart}
\usepackage{amsmath,amsthm,amssymb}
\usepackage{mathtools}
\usepackage[margin=2cm]{geometry}
\usepackage{enumerate}
\usepackage{graphicx}
\usepackage{caption}
\usepackage{xcolor}
\usepackage[hidelinks]{hyperref}
\usepackage{tikz}
\usetikzlibrary{arrows}
\usetikzlibrary[patterns]

\numberwithin{equation}{section}
\usepackage{lineno}
\newcommand{\Bc}{\mathcal{B}}

\newcommand{\cM}{\mathcal{M}}

\renewcommand{\qedsymbol}{$\square$}

\def\opn#1#2{\def#1{\operatorname{#2}}} % to make operators
\opn\chara{char} \opn\length{\ell} \opn\pd{pd} \opn\rk{rk}
\opn\projdim{proj\,dim} \opn\injdim{inj\,dim} \opn\rank{rank}
\opn\depth{depth} \opn\grade{grade} \opn\height{height}
\opn\embdim{emb\,dim} \opn\codim{codim}

\opn\Tr{Tr} \opn\bigrank{big\,rank}
\opn\superheight{superheight}\opn\lcm{lcm}
\opn\trdeg{tr\,deg}%\emph{
	\opn\reg{reg} \opn\lreg{lreg} \opn\ini{in} \opn\lpd{lpd}
	\opn\size{size} \opn\sdepth{sdepth}
	\opn\link{link}\opn\fdepth{fdepth}\opn\lex{lex}\opn\dist{dist}
	\opn\div{div} \opn\Div{Div} \opn\cl{cl} \opn\Cl{Cl}
	\opn\Spec{Spec} \opn\Supp{Supp} \opn\supp{supp} \opn\Sing{Sing}
	\opn\Ass{Ass} \opn\Min{Min}\opn\Mon{Mon}
	\opn\Ann{Ann} \opn\Rad{Rad} \opn\Soc{Soc}
	\opn\Im{Im} \opn\Ker{Ker} \opn\Coker{Coker} \opn\Am{Am}
	\opn\Hom{Hom} \opn\Tor{Tor} \opn\Ext{Ext} \opn\End{End}
	\opn\Aut{Aut} \opn\id{id}
	
	\opn\nat{nat}
	\opn\pff{pf}%   \pf exists already
	\opn\Pf{Pf} \opn\GL{GL} \opn\SL{SL} \opn\mod{mod} \opn\ord{ord}
	\opn\Gin{Gin} \opn\Hilb{Hilb}\opn\sort{sort}
	\opn\aff{aff} \opn
	\con{conv} \opn\relint{relint} \opn\st{st}
	\opn\lk{lk} \opn\cn{cn} \opn\core{core} \opn\vol{vol}
	\opn\link{link} \opn\star{star}\opn\lex{lex}\opn\set{set}
	\opn\gr{gr}
	
	\def\pot#1#2{#1[\kern-0.28ex[#2]\kern-0.28ex]}

	\opn\dirlim{\underrightarrow{\lim}}
	\opn\inivlim{\underleftarrow{\lim}}
	\let\to=\rightarrow
	
	\def\Implies{\ifmmode\Longrightarrow \else
		\unskip${}\Longrightarrow{}$\ignorespaces\fi}
	\def\implies{\ifmmode\Rightarrow \else
		\unskip${}\Rightarrow{}$\ignorespaces\fi}
	\def\iff{\ifmmode\Longleftrightarrow \else
		\unskip${}\Longleftrightarrow{}$\ignorespaces\fi}

	\let\:=\colon
	\let\epsilon\varepsilon
	\let\kappa=\varkappa
	\def\qed{\ifhmode\textqed\fi
		\ifmmode\ifinner\quad\qedsymbol\else\dispqed\fi\fi}
	\def\textqed{\unskip\nobreak\penalty50
		\hskip2em\hbox{}\nobreak\hfil\qedsymbol
		\parfillskip=0pt \finalhyphendemerits=0}
	\def\dispqed{\rlap{\qquad\qedsymbol}}
	
	\opn\dis{dis}
	\def\pnt{{\raise0.5mm\hbox{\large\bf.}}}
	
	\opn\Lex{Lex}
	
        \newtheorem{Theorem}{Theorem}[section]
	\newtheorem{Lemma}[Theorem]{Lemma}
	\newtheorem{Corollary}[Theorem]{Corollary}
	\newtheorem{Proposition}[Theorem]{Proposition}
	\newtheorem{Remark}[Theorem]{Remark}
	
	\newtheorem{Example}[Theorem]{Example}

\begin{document}

        \title[Ideals of $2$-minors of nearly determinantal matrices]{Ideals of $2$-minors of nearly determinantal matrices}

\author[T.~Hibi]{Takayuki Hibi}
\address[Takayuki Hibi]
{Department of Pure and Applied Mathematics, 
Graduate School of Information Science and Technology, 
Osaka University, 
Suita, Osaka 565-0871, Japan}
\email{hibi@math.sci.osaka-u.ac.jp}

\author[F. Navarra]{Francesco Navarra}      
   \address[Francesco Navarra]{Sabanci University, Faculty of Engineering and Natural Sciences, Orta Mahalle, Tuzla 34956, Istanbul, Turkey}	
\email{francesco.navarra@sabanciuniv.edu}

 \author[A. A. Qureshi]{Ayesha Asloob Qureshi}      
    \address[Ayesha Asloob Qureshi]{Sabanci University, Faculty of Engineering and Natural Sciences, Orta Mahalle, Tuzla 34956, Istanbul, Turkey}	
 \email{aqureshi@sabanciuniv.edu, ayesha.asloob@sabanciuniv.edu}

% \author[S.~Saeedi~Madani]{Sara Saeedi Madani}
% \address[Sara Saeedi Madani]
% {Department of Mathematics and Computer Science, Amirkabir University of Technology, Tehran, Iran, and School of Mathematics, Institute for Research in Fundamental Sciences, Tehran, Iran} 
% \email{sarasaeedi@aut.ac.ir, sarasaeedim@gmail.com}

	
	% \subjclass[2010]{13D02, 05E40, 05E45, 05B50, 05C70}
    	  %05C31, graph polynomials, 	13P25  	
    
	\thanks{This research was supported through the program “Oberwolfach Research Fellows” (OWRF) by the Mathematisches Forschungsinstitut Oberwolfach in 2026. The second and third authors are supported by Scientific and Technological Research Council of Turkey T\"UB\.{I}TAK under the Grant 124F113, they are thankful to T\"UB\.{I}TAK for their support. The second author is a member of GNSAGA Indam and he acknowledges their support.} 

\begin{abstract}
Motivated by the study of prime binomial ideals with
Cohen--Macaulay quotients and linear resolutions, we introduce
nearly determinantal matrices, a class of matrices obtained by
repeating variables from the preceding rows in the last row.
We characterize the primality of their ideals of $2$-minors
in terms of an associated directed graph with row labels.
For every prime case, we prove that the quotient is a normal
Cohen--Macaulay domain and determine its graded Betti numbers. Finally, a prime ideal in this
family has a linear resolution if and only if the matrix has
three rows or two columns.
\end{abstract}

        \maketitle

\section*{Introduction}
Determinantal ideals are among the classical and most extensively
studied ideals in commutative algebra. Of particular interest here
are ideals generated by $2$-minors of matrices of variables, since
their generators are quadratic binomials and hence connect
determinantal ideals with the theory of binomial and toric ideals.

Let $X_0$ be a generic $m\times n$ matrix of variables over a field
$K$, where $m,n\ge2$, and let $S$ be the polynomial ring generated
by its entries. It is well known that $I_2(X_0)$ is prime and that
$S/I_2(X_0)$ is Cohen--Macaulay; see, for example,
\cite{BrunsVetter}. Moreover, when $m=2$, the ideal $I_2(X_0)$ has
a linear resolution. Thus, the ideal of $2$-minors of a generic
$2\times n$ matrix simultaneously satisfies three properties:
primality, Cohen--Macaulayness of the quotient, and a linear
resolution. This combination of properties motivated the study,
in \cite{HQS}, of quadratic binomial ideals satisfying the same
three conditions and, in particular, the question of which
matrices of variables give rise to such ideals of $2$-minors.

Accordingly, a matrix $X$ of variables is called a \emph{PCML
matrix} if $I_2(X)$ is prime, $S/I_2(X)$ is Cohen--Macaulay,
and $I_2(X)$ has a linear resolution, where $S$ denotes the
polynomial ring generated by the distinct entries of $X$.
Allowing repeated entries leads to the problem of determining
which patterns of variable identifications produce PCML matrices. The two-row case is closely related to the classical
Kronecker--Weierstrass theory of matrix pencils; see
\cite[Chapter~XII]{Gantmacher}. Over an algebraically closed
field, this theory yields a decomposition of two-row matrices
of linear forms into scroll, Jordan, and nilpotent blocks;
see \cite[Section~3]{CJ}. Concatenations of scroll blocks
define rational normal scrolls, which belong to the classical
classification of varieties of minimal degree
\cite{EH87}. Their defining ideals are prime, their coordinate
rings are Cohen--Macaulay, and their defining ideals have
linear resolutions; see
\cite[Appendix~A2H]{EisenbudSyzygies}.

In Appendix~\ref{app:two-row}, we give a combinatorial treatment
of this classical case for matrices whose entries are variables.
For a two-row matrix $X$ with pairwise distinct entries in each
row and in each column, Theorem~\ref{thm:app-two-row-classification}
gives
\[
X\text{ is PCML}
\quad\Longleftrightarrow\quad
I_2(X)\text{ is prime}
\quad\Longleftrightarrow\quad
\mathcal G_X\text{ contains no directed cycle},
\]
where $\mathcal G_X$ records the repetitions between the two rows.
In the prime case, the path components of $\mathcal G_X$ determine
the scroll blocks.

A natural next question is whether an analogous classification
can be obtained for $m\times n$ matrices with $m>2$. For arbitrary
matrices with repeated variables, however, the possible repetition
patterns become more complicated, and the directed-graph
description of the two-row case does not extend directly.
We therefore consider nearly determinantal matrices, in which
the entries of the first $m-1$ rows are pairwise distinct
variables and the entries of the last row are chosen from these
variables. The repetitions are recorded by a directed graph
whose edges carry row labels. We first classify the patterns
for which the ideal of $2$-minors is prime, and then determine
which of these prime ideals give rise to PCML matrices.

Let $S=K\bigl[x_i^{(j)}:1\leq i\leq n,\ 1\leq j\leq m-1\bigr]$ be a polynomial ring in $(m-1)n$ variables over a field $K$, where
$m,n>1$. We consider matrices of the form
\[
X_{m,n}=
\begin{pmatrix}
x_1^{(1)} & x_2^{(1)} & \cdots & x_n^{(1)}\\
x_1^{(2)} & x_2^{(2)} & \cdots & x_n^{(2)}\\
\vdots & \vdots & & \vdots\\
x_1^{(m-1)} & x_2^{(m-1)} & \cdots & x_n^{(m-1)}\\
y_1 & y_2 & \cdots & y_n
\end{pmatrix},
\]
where, for every $j\in[n]$, the entry $y_j$ is one of the variables
appearing in the first $m-1$ rows; that is,
\[
y_j\in
\bigl\{
x_i^{(k)}:1\leq i\leq n,\ 1\leq k\leq m-1
\bigr\}.
\]
Since our main interest is the primality of $I(X_{m,n})$, we restrict
our attention to matrices satisfying the following conditions:
\begin{enumerate}
    \item[(ND1)] the entries $y_1,\ldots,y_n$ are pairwise distinct;
    \item[(ND2)] no variable is repeated in a column, that is, $y_j\neq x_j^{(i)}$ for all $i\in[m-1]$ and $j\in[n].$
\end{enumerate}
We call such matrices \emph{nearly determinantal matrices}, and denote
by $I(X_{m,n})$ the ideal generated by the $2$-minors of $X_{m,n}$.

These conditions are natural for the study of primality. Indeed, if (ND2) fails, say $y_j=x_j^{(i)}$, then for every $s\neq j$ the $2$-minor determined by the $i$-th row, the last row and the columns $j,s$ is equal to $x_j^{(i)}(y_s-x_s^{(i)})$. Hence $I(X_{m,n})$ is not prime, unless the last row is a copy of the $i$-th row, in which case $I(X_{m,n})$ is the ideal of $2$-minors of a generic $(m-1)\times n$ matrix. Similarly, if (ND1) fails, say $y_p=y_q=z$ with $p\neq q$, then for every $a\in[m-1]$ the $2$-minor determined by the $a$-th row, the last row and the columns $p,q$ is equal to $z(x_p^{(a)}-x_q^{(a)})$, and again $I(X_{m,n})$ is not prime.

The repetitions appearing in the last row admit a natural combinatorial encoding. For every $j\in[n]$, write $y_j=x_{\sigma(j)}^{(\alpha(j))}$
with $\alpha(j)\in[m-1]$ and $\sigma(j)\in[n]$. We associate with $X_{m,n}$ a directed graph $\mathcal G_X$ on the vertex set $[n]$,
having the directed edge $j\longrightarrow \sigma(j)$ for every $j\in[n]$. The value $\alpha(j)$ records the row from which
the repeated variable is taken. A directed cycle
$j_1\longrightarrow j_2\longrightarrow\cdots
\longrightarrow j_r\longrightarrow j_1$ is called \emph{row-distinct} if $\alpha(j_1),\ldots,\alpha(j_r)$ are pairwise distinct.

Our first main result shows that the primality of $I(X_{m,n})$ is
completely determined by the cycle structure of $\mathcal G_X$.

\begin{Theorem}\label{thm:prime-classification}
Let $X_{m,n}$ be a nearly determinantal matrix. Then
$I(X_{m,n})$ is prime if and only if one of the following conditions
holds:
\begin{enumerate}
     \item[{\em (i)}] $m=3$, the integer $n$ is even, and $\mathcal G_X$ is a
    disjoint union of row-distinct directed cycles of length $2$;
    \item[{\em (ii)}] $m\geq4$ and $\mathcal G_X$ consists of a single
    row-distinct directed cycle of length $n$.
\end{enumerate}
\end{Theorem}

The three-row case is exceptional: there exist nearly determinantal
matrices with arbitrarily many columns whose ideals of $2$-minors
are prime. For each fixed $m\ge4$, however, primality forces
$n\le m-1$. In each prime case, we construct a toric parametrization
and a quadratic Gr\"obner basis with squarefree initial monomials.
The corresponding quotient ring is therefore a normal
Cohen--Macaulay domain.

The homological properties of these prime ideals reveal
a further distinction. For $m=3$, write $n=2q$. We identify $I(X_{3,n})$ with the ideal of
$2$-minors of a $2\times3q$ matrix consisting of $q$ scroll blocks of
size three. Thus $I(X_{3,n})$ defines the rational normal scroll of
type $(3,\ldots,3)$, and its graded Betti numbers are given by the
Eagon--Northcott resolution. We also prove that the squarefree
quadratic initial ideal constructed in Section~\ref{subsec:m=3}
has a linear resolution. It follows that the degeneration to this
initial ideal preserves all graded Betti numbers.

For $m\geq 4$, the situation is different. We prove that in this case, the graded Betti numbers of $S/I(X_{m,n})$
coincide with those of the quotient by the ideal of $2$-minors of a
generic $m\times n$ matrix. In particular,
\[
\operatorname{reg} I(X_{m,n})=n
\qquad\text{and}\qquad
\operatorname{pd}_S(S/I(X_{m,n}))=(m-1)(n-1).
\]
Hence, in this case, a prime ideal $I(X_{m,n})$ has a linear resolution
precisely when $n=2$. Combining these results with
Theorem~\ref{thm:prime-classification}, we obtain the following classification of
PCML matrices.
\begin{Theorem}\label{thm:PCML-classification}
Let $X_{m,n}$ be a nearly determinantal matrix. Then $X_{m,n}$ is
PCML if and only if one of the following conditions holds:
\begin{enumerate}
    \item[{\em (i)}] $m=3$, the integer $n$ is even, and $\mathcal G_X$ is a
    disjoint union of row-distinct directed cycles of length $2$;
    \item[{\em (ii)}] $m\geq4$, $n=2$, and $\mathcal G_X$ consists of a single
    row-distinct directed cycle of length $2$.
\end{enumerate}
\end{Theorem}

Equivalently, among the nearly determinantal matrices with prime
ideal of $2$-minors, the PCML property occurs precisely when $m=3$ or
$n=2$. In both cases, we obtain a two-row presentation of the ideal
by scroll blocks: the type is $(3,\ldots,3)$ when $m=3$ and
$(3,1,\ldots,1)$ when $n=2$. These presentations are related to the
path decomposition in Appendix~\ref{app:two-row}.

In Section~\ref{sec:1}, we derive necessary conditions for
primality by analyzing the directed graph $\mathcal G_X$.
In Section~\ref{sec:2}, we prove their sufficiency by constructing
toric parametrizations and squarefree quadratic initial ideals.
For $m=3$, we also identify the resulting ideals with scroll ideals.
Section~\ref{sec:3} determines the graded Betti numbers and
regularity of the prime quotients and proves that the specified
initial ideals in the three-row case have linear resolutions.
Appendix~\ref{app:two-row} gives the two-row classification
in terms of directed paths and rational normal scrolls.

%=%=%=%=%=%=%=%=%=%=%=%=%=%=%=%=%=%=%=%=%=%=%=%=%=%=%=%=%=
%=%=%=%=%=%=%=%=%=%=%=%=%=%=%=%=%=%=%=%=%=%=%=%=%=%=%=%=%=

\section{Directed Cycles and Necessary Conditions for Primality}
\label{sec:1}

We first consider $m=2$. By \textup{(ND1)}, the entries of the last
row are a permutation of the entries of the first row, and this
permutation has no fixed points by \textup{(ND2)}. Thus
$\mathcal G_X$ is a disjoint union of directed cycles of length at
least two. Theorem~\ref{thm:app-two-row-classification} implies that
$I(X_{2,n})$ is not prime. Hence we assume $m\ge3$ throughout this
section.

\begin{Remark}\label{rem:alpha-sigma}
{\em
For every $j\in[n]$, there exist unique indices $\alpha(j)\in[m-1]$ and $
\sigma(j)\in[n]$ such that $y_j=x_{\sigma(j)}^{(\alpha(j))}.$ Thus the last row determines two maps
\[
\alpha:[n]\longrightarrow[m-1]
\qquad\text{and}\qquad
\sigma:[n]\longrightarrow[n].
\]
Condition {\rm (ND2)} is equivalent to
\[
\sigma(j)\neq j
\qquad\text{for every }j\in[n],
\]
whereas Condition {\rm (ND1)} is equivalent to the injectivity of the map
\[
j\longmapsto\bigl(\alpha(j),\sigma(j)\bigr).
\]
}
\end{Remark}

%=%=%=%=%=%=%=%=%=%=%=%=%=%=%=%=%=%=%=%=%=%=%=%=%=%=%=%=%=
%=%=%=%=%=%=%=%=%=%=%=%=%=%=%=%=%=%=%=%=%=%=%=%=%=%=%=%=%=

%=%=%=%=%=%=%=%=%=%=%=%=%=%=%=%=%=%=%=%=%=%=%=%=%=%=%=%=%=

We now derive the restrictions imposed on $\mathcal G_X$ by the
primality of $I(X_{m,n})$. We first show that $\sigma$ must be a permutation, so that $\mathcal G_X$ is a disjoint union of directed cycles. We then prove that these cycles are row-distinct and that any two of them involve the same set of rows. These facts will lead to the two possible configurations occurring in the classification.

%=%=%=%=%=%=%=%=%=%=%=%=%=%=%=%=%=%=%=%=%=%=%=%=%=%=%=%=%=

\begin{Lemma}\label{lem:sigma-bijective}
Let $X_{m,n}$ be a nearly determinantal matrix. If $I(X_{m,n})$ is
prime, then the map $\sigma:[n]\longrightarrow[n]$ is bijective.
\end{Lemma}

\begin{proof}
Since $[n]$ is finite, it is enough to prove that $\sigma$ is
injective. Suppose, on the contrary, that there exist distinct
$p,q\in[n]$ such that $\sigma(p)=\sigma(q)=k.$ Set $a=\alpha(p)$ and $b=\alpha(q).$ Then
\[
y_p=x_k^{(a)}
\qquad\text{and}\qquad
y_q=x_k^{(b)}.
\]
Since the entries $y_1,\ldots,y_n$ are pairwise distinct, we must have
$a\neq b$. Consider the $2$-minor of $X_{m,n}$ determined by the $a$-th row, the
last row, and columns $p,q$. We obtain
\[
x_p^{(a)}y_q-x_q^{(a)}y_p
=
x_p^{(a)}x_k^{(b)}
-
x_q^{(a)}x_k^{(a)}
\in I(X_{m,n}).
\]
On the other hand, the $2$-minor determined by the $a$-th and $b$-th
rows and columns $p,k$ gives
\[
x_p^{(a)}x_k^{(b)}
-
x_k^{(a)}x_p^{(b)}
\in I(X_{m,n}).
\]
Subtracting the latter relation from the former yields $x_k^{(a)}
\bigl(x_p^{(b)}-x_q^{(a)}\bigr)
\in I(X_{m,n}).$ Both factors are nonzero linear forms. Indeed, $x_k^{(a)}$ is a
variable, while $x_p^{(b)}\neq x_q^{(a)}$ since the variables occurring in the first $m-1$ rows are pairwise
distinct. Moreover, $x_k^{(a)}\notin I(X_{m,n})$ and $x_p^{(b)}-x_q^{(a)}\notin I(X_{m,n}),$ because $I(X_{m,n})$ is a homogeneous ideal generated by quadratic
forms and hence contains no nonzero linear form. This contradicts the
primality of $I(X_{m,n})$. Therefore $\sigma$ is injective.
\end{proof}

%=%=%=%=%=%=%=%=%=%=%=%=%=%=%=%=%=%=%=%=%=%=%=%=%=%=%=%=%=

\begin{Remark}\label{rem:sigma-cycles}
{\em 
Suppose that $I(X_{m,n})$ is prime. By Lemma~\ref{lem:sigma-bijective}, the map
\[
\sigma:[n]\longrightarrow[n]
\]
is a permutation. Hence every vertex of $\mathcal G_X$ has indegree and outdegree equal to one, and therefore the connected components of $\mathcal G_X$ are precisely the directed cycles arising from the
cycle decomposition of $\sigma$.

Moreover, by Remark~\ref{rem:alpha-sigma}, $\sigma(j)\neq j$ for every $j\in[n].$ Thus $\sigma$ has no fixed points, and consequently
$\mathcal G_X$ is a disjoint union of directed cycles of length at least $2$. Finally, every row-distinct directed cycle of $\mathcal G_X$ has length at most $m-1$, since the values of $\alpha$ along such a cycle are pairwise distinct elements of $[m-1]$.
}
\end{Remark}

%=%=%=%=%=%=%=%=%=%=%=%=%=%=%=%=%=%=%=%=%=%=%=%=%=%=%=%=%=%=%=%=%=%=%=%=%=%=%=%=%=%=%=%=%=%=%=%=%=%=%=%=%=%=%=%=%=%=%=%=%=%=%=%=%=%=%=%=%=%=%=%=%=%=%=
%=%=%=%=%=%=%=%=%=%=%=%=%=%=%=%=%=%=%=%=%=%=%=%=%=%=%=%=%=%=%=%=%=%=%=%=%=%=%=%=%=%=%=%=%=%=%=%=%=%=%=%=%=%=%=%=%=%=

\begin{Lemma}\label{lem:cycles-row-distinct}
Let $X_{m,n}$ be a nearly determinantal matrix. If $I(X_{m,n})$ is
prime, then every directed cycle of $\mathcal{G}_X$ is row-distinct.
Consequently, every directed cycle of $\mathcal{G}_X$ has length at
most $m-1$.
\end{Lemma}

\begin{proof}
Suppose, to the contrary, that $\mathcal{G}_X$ contains a directed
cycle which is not row-distinct. After relabeling the vertices of this
cycle, we may write it as $1\longrightarrow 2\longrightarrow\cdots\longrightarrow r
\longrightarrow 1.$ Thus, with indices taken modulo $r$, we have $\sigma(j)=j+1$ and hence $y_j=x_{j+1}^{(\alpha(j))}$ for $j=1,\ldots,r$. Since the cycle is not row-distinct, there exist distinct
$p,q\in[r]$ and some $a\in[m-1]$ such that $\alpha(p)=\alpha(q)=a.$ Hence
\[
y_p=x_{p+1}^{(a)}
\qquad\text{and}\qquad
y_q=x_{q+1}^{(a)}.
\]
Let $b\in[m-1]$ with $b\neq a$, and, for simplicity, set
\[
A_j=x_j^{(a)},
\qquad
B_j=x_j^{(b)}.
\]
The $2$-minor determined by the $a$-th row, the last row, and the
columns $p,q$ gives
\[
F:=A_pA_{q+1}-A_qA_{p+1}\in I(X_{m,n}).
\]
Moreover, for every $t$, the $2$-minor determined by the $a$-th and
$b$-th rows and the columns $p,t$ gives
\[
H_t:=A_pB_t-A_tB_p\in I(X_{m,n}).
\]
Set
\[
G:=B_pB_{q+1}-B_qB_{p+1}.
\]
A direct calculation yields
\[
\begin{aligned}
A_p^2G
={}&
A_pB_pH_{q+1}
-A_pB_{p+1}H_q
-A_qB_pH_{p+1}+B_p^2F.
\end{aligned}
\]
Since all the terms on the right-hand side belong to $I(X_{m,n})$, it follows that $A_p^2 G
\in I(X_{m,n}).$ Since $I(X_{m,n})$ is generated in degree $2$, it contains no variables, and since it is prime, one has $A_p^2\notin I(X_{m,n})$. Thus, it is enough to show that we can choose $b$ in such a way that $G\notin I(X_{m,n})$, contradicting the primality of $I(X_{m,n})$.

First suppose that $m\geq4$. Choose $b\in[m-1]\setminus\{a,\alpha(q-1)\}.$ Such a choice is possible since $m-1\geq3$. We claim that the monomial $B_qB_{p+1}$ does not occur as a term of any $2$-minor of $X_{m,n}$. Indeed, a minor involving two of the first $m-1$ rows has, in each term, variables coming from two distinct rows, and therefore cannot contain the product $B_qB_{p+1}$. Thus the only possibility would be a minor involving the $b$-th row and the last row. Such a minor, using columns $s,t$, has the form
$B_s y_t-B_t y_s$, where
$y_t=x_{\sigma(t)}^{(\alpha(t))}$ and the columns
$s,t$ need not lie on the cycle.
Since $\sigma$ is injective, for one of its terms
to be $B_qB_{p+1}$ one must have, up to
interchanging $s$ and $t$, either
\[
t=p,\quad \alpha(p)=b,
\qquad\text{or}\qquad
t=q-1,\quad \alpha(q-1)=b.
\]
 The first is impossible because $\alpha(p)=a\neq b$, and the second is impossible by our choice of $b$. Hence $B_qB_{p+1}$ occurs in no quadratic generator of $I(X_{m,n})$, as claimed. 
Since the degree $2$ component of $I(X_{m,n})$ is
spanned by the $2$-minors, and since
$B_qB_{p+1}\neq B_pB_{q+1}$ because $p\neq q$
and $r\geq 2$, it follows that
$G\notin I(X_{m,n})$.

It remains to consider $m=3$. In this case there are only two possible
values of $\alpha$, say $a$ and $b$. If two consecutive edges of the
cycle come from the same row, say $\alpha(p)=\alpha(p+1)=a,$ we apply the above argument with $q=p+1$. Since $\alpha(q-1)=\alpha(p)=a\neq b,$ the same argument shows that $G\notin I(X_{m,n})$. We are therefore left with the case where no two consecutive edges
come from the same row. Since there are only two rows, the values of
$\alpha$ must alternate around the cycle. As the cycle is not
row-distinct, its length is at least $4$. Choose three consecutive
vertices such that
\[
\alpha(p)=a,\qquad
\alpha(p+1)=b,\qquad
\alpha(p+2)=a,
\]
and set $q=p+2$. Then $G=B_pB_{p+3}-B_{p+2}B_{p+1}.$ We claim that the monomial $B_{p+1}B_{p+2}$ does not occur in any $2$-minor. As above, it cannot occur in a minor of the first two rows. If it occurred in a minor between
the $b$-th row and the last row, then one of the two factors would have
to be an entry $y_t=x_{t+1}^{(\alpha(t))}$ of the last row. Since $\sigma$ is injective, there are only two
possibilities. Taking $t=p$ is impossible since
$\alpha(p)=a\neq b$, while taking $t=p+1$ would require using the same
column twice and hence does not define a $2$-minor. Thus
$B_{p+1}B_{p+2}$ occurs in no quadratic generator, and again $G\notin I(X_{m,n}).$
\end{proof}

%=%=%=%=%=%=%=%=%=%=%=%=%=%=%=%=%=%=%=%=%=%=%=%=%=%=%=%=%=
%=%=%=%=%=%=%=%=%=%=%=%=%=%=%=%=%=%=%=%=%=%=%=%=%=%=%=%=%=

\begin{Lemma}\label{lem:cycles-same-rows}
Let $X_{m,n}$ be a nearly determinantal matrix and suppose that
$I(X_{m,n})$ is prime. Then any two directed cycles of
$\mathcal{G}_X$ involve the same set of rows. More precisely, if
$C$ and $D$ are directed cycles of $\mathcal{G}_X$, then
\[
\{\alpha(j):j\in V(C)\}
=
\{\alpha(j):j\in V(D)\}.
\]
Consequently, all directed cycles of $\mathcal{G}_X$ have the same
length. In particular, if their common length is $r$, then
$r\leq m-1$ and $r$ divides $n$.
\end{Lemma}

\begin{proof}
Set $I=I(X_{m,n})$. By Remark~\ref{rem:sigma-cycles},
$\mathcal G_X$ is a disjoint union of directed cycles of length
at least $2$. Suppose, contrary to the assertion, that there exist
two directed cycles having different row sets. After simultaneously
relabeling the columns and the subscripts of the variables, we may
assume that these cycles are
\[
C:
1\longrightarrow 2\longrightarrow\cdots\longrightarrow r
\longrightarrow 1
\quad \text{and}\quad
D:
r+1\longrightarrow r+2\longrightarrow\cdots
\longrightarrow r+s\longrightarrow r+1.
\]
For $i\in[r]$, define
\[
i^+=
\begin{cases}
i+1,&1\leq i<r,\\
1,&i=r.
\end{cases}
\]
Then $\sigma(i)=i^+$, and
\begin{equation}\label{eq:1}
y_i=x_{i^+}^{(\alpha(i))}
\qquad\text{for every }i\in[r].
\end{equation}
Set $\mathcal R(C)=\{\alpha(1),\ldots,\alpha(r)\}$ and
$\mathcal R(D)=\{\alpha(r+1),\ldots,\alpha(r+s)\}$.
By assumption, $\mathcal R(C)\neq\mathcal R(D)$.
After interchanging $C$ and $D$, if necessary, we may choose
\begin{equation}\label{eq:2}
a\in\mathcal R(C)\setminus\mathcal R(D).
\end{equation}
We claim that
\begin{equation}\label{eq:3}
x_i^{(a)}y_{r+1}
-
x_{i^+}^{(a)}x_{r+1}^{(\alpha(i))}
\in I
\qquad\text{for every }i\in[r].
\end{equation}
Indeed, the $2$-minor determined by the $a$-th row, the last row,
and columns $i,r+1$ gives
$x_i^{(a)}y_{r+1}-x_{r+1}^{(a)}y_i\in I$.
Using~\eqref{eq:1}, this becomes
\begin{equation}\label{eq:4}
x_i^{(a)}y_{r+1}
-
x_{r+1}^{(a)}x_{i^+}^{(\alpha(i))}
\in I.
\end{equation}
If $\alpha(i)=a$, then~\eqref{eq:4} is precisely~\eqref{eq:3}.
If $\alpha(i)\neq a$, the $2$-minor determined by the $a$-th and
$\alpha(i)$-th rows and columns $r+1,i^+$ gives
\begin{equation}\label{eq:5}
x_{r+1}^{(a)}x_{i^+}^{(\alpha(i))}
-
x_{i^+}^{(a)}x_{r+1}^{(\alpha(i))}
\in I.
\end{equation}
Adding~\eqref{eq:4} and~\eqref{eq:5} yields~\eqref{eq:3}. Set
\[
M_C=\prod_{i=1}^r x_i^{(a)}
\qquad\text{and}\qquad
N_C=\prod_{i=1}^r x_{r+1}^{(\alpha(i))}.
\]
Relation~\eqref{eq:3} gives, in the quotient ring $S/I$,
\[
x_i^{(a)}y_{r+1}
=
x_{i^+}^{(a)}x_{r+1}^{(\alpha(i))}
\qquad\text{for every }i\in[r].
\]
Multiplying these equalities for $i=1,\ldots,r$, we obtain
\[
\left(\prod_{i=1}^r x_i^{(a)}\right)y_{r+1}^r
=
\left(\prod_{i=1}^r x_{i^+}^{(a)}\right)
\left(\prod_{i=1}^r x_{r+1}^{(\alpha(i))}\right)
\quad\text{in }S/I.
\]
Since $i\mapsto i^+$ permutes the vertices of $C$,
\[
\prod_{i=1}^r x_{i^+}^{(a)}
=
\prod_{i=1}^r x_i^{(a)}
=
M_C.
\]
Consequently,
\begin{equation}\label{eq:6}
M_C\bigl(y_{r+1}^r-N_C\bigr)\in I.
\end{equation}
Since $I$ is generated in degree $2$, it contains no variables,
and since it is prime, one has $M_C\notin I$.
Next, we prove that $y_{r+1}^r-N_C\notin I$, which contradicts
the primality of $I$. Define a $K$-algebra homomorphism
$\varphi:S\longrightarrow K$ by
\begin{equation}\label{eq:7}
\varphi\bigl(x_j^{(b)}\bigr)=
\begin{cases}
1,
&
r+1\leq j\leq r+s
\text{ and }b\neq a,
\\[1mm]
0,
&
\text{otherwise}.
\end{cases}
\end{equation}
We claim that $I\subseteq\Ker\varphi$.
Indeed, if $j\in V(D)=\{r+1,\ldots,r+s\}$, then
$\sigma(j)\in V(D)$. Moreover, by~\eqref{eq:2},
$\alpha(j)\neq a$. Therefore
$\varphi(y_j)
=\varphi\left(x_{\sigma(j)}^{(\alpha(j))}\right)=1$.
Thus every column indexed by a vertex of $D$ is mapped to the same
column vector
\[
\begin{pmatrix}
1\\
\vdots\\
1\\
0\\
1\\
\vdots\\
1
\end{pmatrix},
\]
where the zero occurs in the $a$-th position among the first $m-1$
entries. On the other hand, if $j\notin V(D)$, then
$\sigma(j)\notin V(D)$, because $D$ is a cycle of the permutation
$\sigma$. Hence
$\varphi(y_j)
=\varphi\left(x_{\sigma(j)}^{(\alpha(j))}\right)=0$.
Therefore every column outside $D$ is mapped to the zero column.
Hence any two columns of $\varphi(X_{m,n})$ are either equal or
one of them is zero. Consequently, $I\subseteq\Ker\varphi$.
Moreover, since $r+1\in V(D)$ and $a\notin\mathcal R(D)$,
we have $\varphi(y_{r+1})=1$.
On the other hand, $a\in\mathcal R(C)$, so there exists some
$i\in[r]$ for which $\alpha(i)=a$. Therefore $N_C$ is divisible by
$x_{r+1}^{(a)}$, and hence $\varphi(N_C)=0$.
Consequently, $\varphi\bigl(y_{r+1}^r-N_C\bigr)=1$, and thus
$y_{r+1}^r-N_C\notin I$.
Relation~\eqref{eq:6} is therefore a product of two elements not
belonging to $I$, contradicting the primality of $I$.
We conclude that any two directed cycles of $\mathcal G_X$
involve the same set of rows.

Finally, by Lemma~\ref{lem:cycles-row-distinct}, every directed cycle
is row-distinct. Hence the length of a cycle is equal to the
cardinality of its row set. Since all directed cycles have the same
row set, they have a common length, say $r$, with $r\leq m-1$.
Since the directed cycles form a partition of the vertex set $[n]$,
their common length $r$ divides $n$.
\end{proof}

%=%=%=%=%=%=%=%=%=%=%=%=%=%=%=%=%=%=%=%=%=%=%=%=%=%=%=%=%=
%=%=%=%=%=%=%=%=%=%=%=%=%=%=%=%=%=%=%=%=%=%=%=%=%=%=%=%=%=%=%=%=%=%=%=%=%=%=%=%=%=%=%=%=%=%=%=%=%=

\begin{Proposition}\label{prop:m-at-least-four-necessary}
Let $m\geq4$. If $I(X_{m,n})$ is prime, then
$\mathcal G_X$ consists of a single row-distinct directed cycle of
length $n$. In particular, if $4\leq m\leq n$, then
$I(X_{m,n})$ is not prime.
\end{Proposition}

\begin{proof}
Set $I=I(X_{m,n})$. By Remark~\ref{rem:sigma-cycles},
$\mathcal G_X$ is a disjoint union of directed cycles of length
at least $2$. Suppose, contrary to the assertion, that
$\mathcal G_X$ contains at least two directed cycles.
By Lemma~\ref{lem:cycles-row-distinct}, every directed cycle
is row-distinct, and by Lemma~\ref{lem:cycles-same-rows}, all directed
cycles involve the same set of rows. In particular, all of them have
the same length, say $r$. Choose two such cycles. Since they involve
the same set of rows, we may choose a row $a\in[m-1]$ which occurs
on both cycles. After simultaneously relabeling the columns and
their subscripts, and cyclically choosing the initial vertex of
each cycle, we may assume that the two cycles are
\[
C:
1\longrightarrow2\longrightarrow\cdots\longrightarrow r
\longrightarrow1
\quad \text{and}\quad
D:
r+1\longrightarrow r+2\longrightarrow\cdots
\longrightarrow2r\longrightarrow r+1,
\]
and that $\alpha(1)=\alpha(r+1)=a$.
It follows that $\sigma(1)=2$, $\sigma(r+1)=r+2$. Hence
\begin{equation}\label{eq:8}
y_1=x_2^{(a)},
\qquad
y_{r+1}=x_{r+2}^{(a)}.
\end{equation}
Since $m\geq4$, there are at least three rows among the first $m-1$
rows. We may therefore choose
\begin{equation}\label{eq:9}
b\in[m-1]\setminus\{a,\alpha(2r)\}.
\end{equation}
%Here $2r$ is the predecessor of $r+1$ on the cycle $D$, so that \[\sigma(2r)=r+1.\]

We now construct a relation in $I$ that contradicts its primality.
The $2$-minor determined by the $a$-th and $b$-th rows and columns
$1,2$ gives
\begin{equation}\label{eq:10}
f_1=
x_1^{(a)}x_2^{(b)}
-
x_2^{(a)}x_1^{(b)}
\in I.
\end{equation}
The $2$-minor determined by the $b$-th row, the last row, and columns
$1,r+1$ gives, by~\eqref{eq:8},
\begin{equation}\label{eq:11}
f_2
=
x_1^{(b)}y_{r+1}
-
x_{r+1}^{(b)}y_1
=
x_1^{(b)}x_{r+2}^{(a)}
-
x_{r+1}^{(b)}x_2^{(a)}
\in I.
\end{equation}
Finally, the $2$-minor determined by the $a$-th and $b$-th rows and
columns $1,r+2$ gives
\begin{equation}\label{eq:12}
f_3=
x_1^{(a)}x_{r+2}^{(b)}
-
x_{r+2}^{(a)}x_1^{(b)}
\in I.
\end{equation}
A direct computation using~\eqref{eq:10}--\eqref{eq:12} yields
\begin{equation}\label{eq:13}
x_{r+1}^{(b)}f_1
-
x_1^{(b)}f_2
-
x_1^{(b)}f_3
=
x_1^{(a)}
\left(
x_2^{(b)}x_{r+1}^{(b)}
-
x_1^{(b)}x_{r+2}^{(b)}
\right)
\in I.
\end{equation}
Since $I$ is generated in degree $2$, one has $x_1^{(a)}\notin I$.
We claim that the second factor in~\eqref{eq:13} does not belong to
$I$ either. Set
\[
Q=
x_2^{(b)}x_{r+1}^{(b)}
-
x_1^{(b)}x_{r+2}^{(b)}.
\]
Define a $K$-algebra homomorphism $\varphi:S\longrightarrow K$ by
\[
\varphi\bigl(x_2^{(b)}\bigr)
=
\varphi\bigl(x_{r+1}^{(b)}\bigr)=1
\]
and by sending every other variable to zero.
We first verify that every entry in the last row is mapped to zero.
Since $\sigma$ is a permutation, the only last-row entry which could
equal $x_2^{(b)}$ is $y_1=x_2^{(\alpha(1))}$.
However, $\alpha(1)=a\neq b$, and hence $y_1\neq x_2^{(b)}$.
Similarly, the only last-row entry which could equal
$x_{r+1}^{(b)}$ is $y_{2r}=x_{r+1}^{(\alpha(2r))}$,
because $\sigma(2r)=r+1$.
By the choice of $b$ in~\eqref{eq:9}, $\alpha(2r)\neq b$,
and therefore $y_{2r}\neq x_{r+1}^{(b)}$.
It follows that $\varphi(y_j)=0$ for every $j\in[n]$.
Thus, in the matrix $\varphi(X_{m,n})$, the only nonzero entries are
the two entries in the $b$-th row and columns $2$ and $r+1$.
Consequently, all $2$-minors of $\varphi(X_{m,n})$ vanish,
and hence $I\subseteq\Ker\varphi$.
On the other hand, $\varphi(Q)=1$. Therefore $Q\notin I$.

Relation~\eqref{eq:13} is thus a product of two elements not belonging
to $I$, contradicting the primality of $I$.
Hence $\mathcal G_X$ contains only one directed cycle.
Since $\sigma$ is a permutation of $[n]$, this unique cycle contains
all $n$ vertices and therefore has length $n$.
By Lemma~\ref{lem:cycles-row-distinct}, it is row-distinct.
Hence its $n$ edges involve $n$ distinct rows among the first $m-1$
rows, so $n\leq m-1$. Therefore $m>n$.
\end{proof}

%=%=%=%=%=%=%=%=%=%=%=%=%=%=%=%=%=%=%=%=%=%=%=%=%=%=%=%=%=%=%=%=%=%=%=%=%=%=%=%=%=%=%=%=%=%=%=%=%=%=%=%=%=%=%=%=%=%=%=%=%=%=%=%=%=%=%=%=%=%=%=%=%=%=%=%=%=%=%=%=%=%=%=%=
\begin{Proposition}\label{prop:m-equals-three-necessary}
Suppose that $m=3$. If $I(X_{3,n})$ is prime, then $n$ is even and
$\mathcal G_X$ is a disjoint union of row-distinct directed cycles of
length $2$.
\end{Proposition}

\begin{proof}
Suppose that $I(X_{3,n})$ is prime. By Remark~\ref{rem:sigma-cycles},  $\mathcal G_X$ is a disjoint union of directed
cycles of length at least 2. By Lemma~\ref{lem:cycles-row-distinct}, every directed cycle of
$\mathcal G_X$ is row-distinct. Since $m=3$, there are only two
available rows, namely the rows indexed by $1$ and $2$. Therefore a
row-distinct directed cycle has length exactly $2$. In particular, the cycles
partition the vertex set $[n]$ into pairs, and hence $n$
is even.
\end{proof}

%=%=%=%=%=%=%=%=%=%=%=%=%=%=%=%=%=%=%=%=%=%=%=%=%=%=%=%=%=%=%=%=%=%=%=%=%=%=%=%=%=%=%=%=%=%=%=%=%=%=%=%=%=%=%=%=%=%=%=

As an immediate consequence of Proposition~\ref{prop:m-at-least-four-necessary} and Proposition~\ref{prop:m-equals-three-necessary}, we obtain the following result.

\begin{Corollary}\label{cor:main-classification}
Let $X_{m,n}$ be a nearly determinantal matrix. If $I(X_{m,n})$ is prime, then one of the following conditions
holds:
\begin{enumerate}
    \item[{\em(i)}] $m=3$, $n$ is even, and $\mathcal G_X$ is a
    disjoint union of row-distinct directed cycles of length $2$;

    \item[{\em(ii)}] $m\geq4$ and $\mathcal G_X$ consists of a single
    row-distinct directed cycle of length $n$.
\end{enumerate}
    
\end{Corollary}

%=%=%=%=%=%=%=%=%=%=%=%=%=%=%=%=%=%=%=%=%=%=%=%=%=%=%=%=%=%=%=%=%=%=%=%=%=%=%=%=%=%=%=
%=%=%=%=%=%=%=%=%=%=%=%=%=%=%=%=%=%=%=%=%=%=%=%=%=%=%=%=%=%=%=%=%=%=%=%=%=%=%=%=%=%=%=%=%=%=%=%=%=%=%=%=%=%=%=%=%=%=%=%=%=%=%=

\section{Primality and Cohen–Macaulayness via Gröbner Bases}\label{sec:2}

In the previous section, we obtained necessary conditions
for the primality of $I(X_{m,n})$, as summarized in
Corollary~1.8. We now prove that these conditions are also
sufficient. In each case, we construct a toric parametrization
and use Gr\"obner bases to establish the primality of
$I(X_{m,n})$ and the Cohen--Macaulayness of $S/I(X_{m,n})$.

\subsection{The case $m=3$}\label{subsec:m=3}

By Corollary~1.8, we need to consider the case where the number
of columns is even and $\mathcal G_X$ is a disjoint union of row-distinct
directed cycles of length $2$. We write the number of columns
as $2n$. After a suitable relabeling of the columns and
variables, the corresponding matrix has the form
\begin{equation}\label{matrixX}
    X=
\begin{pmatrix}
x_1 & \cdots & x_n & x_{n+1} & \cdots & x_{2n}\\
y_1 & \cdots & y_n & y_{n+1} & \cdots & y_{2n}\\
y_{n+1} & \cdots & y_{2n} & x_1 & \cdots & x_n
\end{pmatrix}.
\end{equation}
Let $S=K[x_1,\ldots,x_{2n},y_1,\ldots,y_{2n}]$. Consider the
$2\times3n$ matrix
\[
H=(H_1\mid\cdots\mid H_n),
\qquad
H_i=
\begin{pmatrix}
x_{n+i}&x_i&y_{n+i}\\
x_i&y_{n+i}&y_i
\end{pmatrix},
\qquad 1\le i\le n.
\]
We show that $I(X)=I_2(H)$. Set $a_i=x_{n+i}, b_i=x_i, c_i=y_{n+i}, d_i=y_i.$ Then
\[
X=
\begin{pmatrix}
b_1&\cdots&b_n&a_1&\cdots&a_n\\
d_1&\cdots&d_n&c_1&\cdots&c_n\\
c_1&\cdots&c_n&b_1&\cdots&b_n
\end{pmatrix}.
\]
The ideal $I_2(H)$ is generated by the following six types of binomials:
\[
\begin{gathered}
a_i b_j-a_j b_i,\qquad
a_i c_j-b_i b_j,\qquad
a_i d_j-b_i c_j,\\
b_i c_j-b_j c_i,\qquad
b_i d_j-c_i c_j,\qquad
c_i d_j-c_j d_i,
\end{gathered}
\qquad 1\le i,j\le n.
\]
Apart from the third type, each is, up to sign, a $2$-minor of $X$.
For the third type, one has
\[
a_i d_j-b_i c_j
=(a_i d_j-b_j c_i)-(b_i c_j-b_j c_i),
\]
and both terms on the right belong to $I(X)$. Thus
$I_2(H)\subseteq I(X)$.

Conversely, every $2$-minor of $X$ is, up to sign, either one of the
six types above or one of the binomials on the left-hand sides of
the following identities:
\[
\begin{aligned}
b_i d_j-b_j d_i
  &=(b_i d_j-c_i c_j)-(b_j d_i-c_i c_j),\\
a_i c_j-a_j c_i
  &=(a_i c_j-b_i b_j)-(a_j c_i-b_i b_j),\\
a_i d_j-b_j c_i
  &=(a_i d_j-b_i c_j)+(b_i c_j-b_j c_i).
\end{aligned}
\]
The right-hand sides belong to $I_2(H)$, so $I(X)\subseteq I_2(H)$.
Hence $I(X)=I_2(H)$.

Each block $H_i$ is a scroll block of size three, in the
terminology of Appendix~\ref{app:two-row}, and different
blocks have disjoint sets of variables. Equivalently, the graph
$\mathcal G_H$ of Appendix~\ref{app:two-row} consists of $n$ directed
paths, each on three vertices. Thus $I(X)=I_2(H)$ defines the rational normal scroll of type
$(3,\ldots,3)$. By Theorem~\ref{thm:app-two-row-classification},
$S/I(X)$ is a Cohen--Macaulay domain and $I(X)$ has a
$2$-linear resolution.

We next prove that the $2$-minors of the original matrix $X$ form a
Gr\"obner basis for the lexicographic order specified below. Their
initial monomials are squarefree. In
Section~\ref{subec:3.2}, we prove that this initial
ideal itself has a linear resolution and that its quotient has the
same graded Betti numbers as $S/I(X)$.

%=%=%=%=%=%=%=%=%=%=%=%=%=%=%=%=%=%=%=%=%=%=%=%=%=%=%=%=

Let $T=K[a^{\pm1},b^{\pm1},x_2,\ldots,x_n]$ be the Laurent polynomial ring in $n+1$ variables over $K$.  We introduce the $K$-algebra homomorphism $\pi : S \to T$ by setting
\begin{eqnarray*}
& \pi(x_1) = a, \quad \pi(x_{n+1}) = a^2/b, \quad \pi(y_1) = b^2/a, \quad \pi(y_{n+1}) = b, &\\
& \pi(x_{i}) = x_i, \quad \pi(x_{n+i}) = (a/b) x_i, \quad 2 \leq i \leq n,& \\
& \pi(y_{i}) = (b/a)^2 x_i, \quad \pi(y_{n+i}) = (b/a)x_i, \quad 2 \leq i \leq n.&
\end{eqnarray*}
Let $\Ker({\pi}) \subset S$ denote the kernel of $\pi$.

    \begin{Lemma}
        One has $I(X) \subseteq \Ker(\pi)$.
    \end{Lemma}

    Let $<_{\lex}$ denote the lexicographic order \cite[Example 2.1.2 (a)]{HHgtm260} on $S$ induced by the ordering
    \[
    x_{n+1} > \cdots > x_{2n} > y_1 > \cdots > y_n > x_1 > \cdots > x_n > y_{n+1} > \cdots > y_{2n} 
    \]
of variables.  Let $\Omega$ denote the set of $2$-minors of $X$ and set
\[
{\rm in}_{<_{\lex}}(\Omega) := \{{\rm in}_{<_{\lex}}(u) : u \in \Omega\}.
\]

    \begin{Lemma}
    Each monomial belonging to ${\rm in}_{<_{\lex}}(\Omega)$ is squarefree.
    \end{Lemma}
    
A monomial $u=\prod_{i=1}^{2n} x_i^{p_i} \prod_{i=1}^{2n} y_i^{q_i}$ of $S$ is called {\em standard} if $u$ is not divisible by any monomial in ${\rm in}_{<_{\lex}}(\Omega)$.  Let $\cM$ denote the set of standard monomials of $S$.

\begin{Lemma}
    If $u \in \cM, v \in \cM$ and $u \neq v$, then $\pi(u) \neq \pi(v)$.
\end{Lemma}

\begin{proof}
Let $\supp(u)$ denote the support of $u \in \cM$.  Suppose that
\begin{itemize}
    \item 
    $\supp(u) \cap \supp(v) = \emptyset$,
    \item 
    $\deg_{\,x_i}{\pi(u)}=\deg_{\,x_i}{\pi(v)}$ for each $2 \leq i \leq n$,
    \item
    $\deg{u}=\deg{v}$.
\end{itemize}

\smallskip

\noindent
{\bf (First Step)} Let $x_{n+1} \in \supp(u)$.  It then follows that
\[
\supp(u) \subset \{x_1, x_{n+1}, \ldots, x_{2n} \}.
\]
One has $\deg_{\,a}{\pi(u)} > \deg{u}$.  Since 
% $\deg_{\,x_i}{\pi(u)}=\deg_{\,x_i}{\pi(v)}$ for each $2 \leq i \leq n$, 
$x_{n+1} \not\in \supp(v)$,
it follows that $\deg_{\,a}{\pi(v)} \leq \deg{v}$.  Since $\deg{u} = \deg{v}$, one has $\deg_{\,a}{\pi(u)} > \deg_{\,a}{\pi(v)}$.  Thus $\pi(u) \neq \pi(v)$.

Let $1 < i \leq n$.  Suppose that $x_{n+i} \in \supp(u)$ and that none of $x_{n+1}, x_{n+2}, \ldots, x_{n+i-1}$ belongs to $\supp(u) \cup \supp(v)$.  It then follows that  
\[
\supp(u) \subset \{x_1, \ldots, x_{i}, x_{n+i}, \ldots, x_{2n} \}.
\]
Since none of $x_{n+1}, x_{n+2}, \ldots, x_{n+i-1}, x_{n+i}$ belongs to $\supp(v)$, and since $\deg_{\,x_i}{\pi(u)}=\deg_{\,x_i}{\pi(v)}$ for each $2 \leq i \leq n$ and $\deg{u} = \deg{v}$, one has $\deg_{\,a}{\pi(u)} > \deg_{\,a}{\pi(v)}$.  Thus $\pi(u) \neq \pi(v)$.

\medskip

\noindent
{\bf (Second Step)}  Suppose that none of $x_{n+1}, x_{n+2}, \ldots, x_{2n}$ belongs to $\supp(u) \cup \supp(v)$.  

Let $y_1 \in \supp(u)$.  It then follows that  
\[
\supp(u) \subset \{y_1, \ldots, y_n, y_{n+1}\}.
\] 
One has $\deg_{\,b}{\pi(u)} > \deg_{\,b}{\pi(v)}$.  Thus $\pi(u) \neq \pi(v)$.

Let $1 \leq i \leq n$.  Suppose that $y_i \in \supp(u)$ and that none of $y_1, y_2, \ldots, y_{i-1}$ belongs to $\supp(u) \cup \supp(v)$.  It then follows that  
\[
\supp(u) \subset \{y_i, \ldots, y_n, y_{n+1}, \ldots, y_{n+i}\}.
\]  
One has $\deg_{\,b}{\pi(u)} > \deg_{\,b}{\pi(v)}$.  Thus $\pi(u) \neq \pi(v)$.

\medskip

\noindent
{\bf (Third Step)}  Suppose that none of $x_{n+1}, x_{n+2}, \ldots, x_{2n}, y_1, \ldots, y_n$ belongs to $\supp(u) \cup \supp(v)$.  

Let $x_1 \in \supp(u)$.  It then follows that  
\[
\supp(u) \subset \{x_1, \ldots, x_n, y_{n+1}\}.
\] 
One has $\deg_{\,a}{\pi(u)} > \deg_{\,a}{\pi(v)}$.  Thus $\pi(u) \neq \pi(v)$.

Let $1 \leq i \leq n$.  Suppose that $x_i \in \supp(u)$ and that none of $x_1, x_2, \ldots, x_{i-1}$ belongs to $\supp(u) \cup \supp(v)$.  It then follows that  
\[
\supp(u) \subset \{x_i, \ldots, x_n, y_{n+1}, \ldots, y_{n+i}\}.
\]
It follows that $y_{n+i} \in \supp(v)$ and that none of $y_{n+1}, \ldots, y_{n+i-1}$ belongs to $\supp(u)$.  
One has $\deg_{\,b}{\pi(u)} < \deg_{\,b}{\pi(v)}$.  Thus $\pi(u) \neq \pi(v)$.

\medskip

\noindent
{\bf (Fourth Step)}  Let $\supp(u) \cup \supp(v) \subset \{y_{n+1}, \ldots, y_{2n}\}$.  Since $y_{n+1}, \ldots, y_{2n}$ are algebraically independent, the desired result follows. 
\hspace{9.1cm}
\end{proof}

A standard fact on Gr\"obner bases
\cite[Proposition~2.2.5]{HHgtm260},
together with the normality criterion
\cite[Proposition~13.15]{Stu}
and Hochster's theorem
\cite[Theorem~1]{Hochster1972},
guarantees that

\begin{Corollary}
    The set of $2$-minors of $X$ is a Gr\"obner basis of $\Ker(\pi)$ with respect to $<_{\lex}$ and its initial monomials are
squarefree. In particular, $I(X)=\Ker(\pi)$
and $S/I(X)$ is a normal and Cohen--Macaulay domain.
\end{Corollary}

This completes the proof of Theorem~\ref{thm:prime-classification}(i).

%=%=%=%=%=%=%=%=%=%=%=%=%=%=%=%=%=%=%=%=%=%=%=%=%=%=%=%=%=%=

\subsection{The case $m\geq 4$}\label{sec:m=4} We now turn to the discussion of the ideal of $2$-minors $I(X)$ of an $m\times n$ nearly determinantal matrix such that
$\mathcal G_X$ consists of a single row-distinct directed
cycle of length $n$, so that $m>n$. After a relabeling
of the columns and of the first $m-1$ rows, we may
assume that
\[
X=
\begin{pmatrix}
x_{1}^{(1)} & x_{2}^{(1)} & \cdots & x_{n-1}^{(1)} & x_{n}^{(1)} \\
\vdots & \vdots &       & \vdots & \vdots\\
x_{1}^{(n)} & x_{2}^{(n)} & \cdots & x_{n-1}^{(n)} & x_{n}^{(n)}\\
x_{1}^{(n+1)} & x_{2}^{(n+1)} & \cdots & x_{n-1}^{(n+1)} & x_{n}^{(n+1)} \\
\vdots & \vdots &       & \vdots & \vdots\\
x_{1}^{(m-1)} & x_{2}^{(m-1)} & \cdots & x_{n-1}^{(m-1)} & x_{n}^{(m-1)}\\
x_2^{(2)} & x_3^{(3)} & \cdots & x_{n}^{(n)} & x_1^{(1)}
\end{pmatrix}.
\]
Let $S=K[\{x_j^{(r)}\}_{1 \leq r\leq m-1, \, 1\leq j \leq n}]$ denote the polynomial ring in $(m-1)n$ variables over a field $K$ and $T=K[a_1^{\pm 1}, \ldots, a_n^{\pm 1}, t_{n+1}, \ldots, t_{m-1}]$ the Laurent polynomial ring in $m-1$ variables over $K$.  We introduce the map $\pi : S \to T$ by setting
\begin{eqnarray}
\label{abcdefg}   
\pi(x_j^{(i)})= \frac{\,a_{j+1}\,}{a_{i+1}} \, a_{i}, \quad 1 \leq i \leq n,
\qquad \pi(x_j^{(r)}) = t_ra_{j+1}, \quad n < r < m,
\end{eqnarray}
where $1 \leq j \leq n$ and $a_{n+1} = a_1$.  Let $\Ker(\pi) \subset S$ denote the kernel of $\pi$. 

\begin{Lemma}
    One has $I(X) \subseteq \Ker(\pi)$.
\end{Lemma}

Let $<_{\lex}$ denote the lexicographic order on $S$ induced by the ordering

\begin{eqnarray*}
x_{1}^{(1)}< x_{2}^{(2)} < \cdots < x_{n}^{(n)} & < & x_1^{(n+1)} < x_2^{(n+1)} < \cdots < x_n^{(n+1)} \\
& < & \cdots
\\
& < & x_1^{(m-1)} < x_2^{(m-1)} < \cdots < x_n^{(m-1)} 
\\ & < &
x_2^{(1)} < x_3^{(1)} < \cdots < x_n^{(1)} \\ & < & 
x_1^{(2)} < x_3^{(2)}  < \cdots < x_n^{(2)} \\ & < & 
x_1^{(3)} < \cdots < x_{n-2}^{(n-1)} < x_n^{(n-1)} \\ & < & x_1^{(n)} < \cdots < x_{n-1}^{(n)}
\end{eqnarray*}
of variables.  For $i\in[n]$, we use cyclic notation
\[
i^+=
\begin{cases}
i+1,& i<n,\\
1,& i=n,
\end{cases}
\qquad
i^-=
\begin{cases}
i-1,& i>1,\\
n,& i=1.
\end{cases}
\]

\noindent
Let $\Omega$ denote the set of all $2$-minors of $X$.  We enlarge
$\Omega$ by adding the following quadratic binomials.  For
$p,q\in[n]$ and $r \in [m-1]$, set
\[
F_{pqr}
:=
x_q^{(p)}x_p^{(r)}
-
x_{q^+}^{(q^+)}x_{p^-}^{(r)}.
\]
Let
\[
\mathcal A
=
\left\{
F_{pqr}:
p,q\in[n], r \in [m-1],\
r\neq p,\
q\notin\{p,p^-\}
\right\},
\]
and let $\Bc =\Omega\cup\mathcal A$. Each $F_{pqr}$ belongs to $I(X)$. Indeed,
\[
F_{pqr}
=
\left(
x_q^{(p)}x_p^{(r)}
-
x_p^{(p)}x_q^{(r)}
\right) +
\left(
x_p^{(p)}x_q^{(r)}
-
x_{q^+}^{(q^+)}x_{p^-}^{(r)}
\right),
\]
and each summand on the right-hand side is, up to sign,
a $2$-minor of $X$ or zero.
The restrictions in the definition of $\mathcal A$ remove precisely
the cases in which $F_{pqr}$ is already a $2$-minor (or is zero).
In fact,
\begin{itemize}
\item if $q=p^-$, then $F_{p,p^-,r}
=
x_{p^-}^{(p)}x_p^{(r)}
-
x_p^{(p)}x_{p^-}^{(r)}$ is a $2$-minor of rows $p,r$;

\item if $q=p$, then $F_{p,p,r}
=
x_p^{(p)}x_p^{(r)}
-
x_{p^+}^{(p^+)}x_{p^-}^{(r)}$ is, up to sign, a $2$-minor involving row $r$ and the last row;

\item if $r=p$, then $F_{pqp}$ is either zero or, up to sign,
a $2$-minor involving row $p$ and the last row.
\end{itemize}
Thus $\mathcal A$ consists exactly of the additional quadratic
binomials in $\Bc$ which do not already belong to $\Omega$.
One has $|\mathcal A|=n(n-2)(m-2).$

%=%=%=%=%=%=%=%=%=%=%=%=%=%=%=%=%=%=%=%=%=%=%=%=
\begin{Example}\label{example}
{\em 
Let $m=8$ and $n=5$. The additional binomials $F_{pqr}$ are indexed by
\[
p\in[5],\qquad
q\in[5]\setminus\{p,p^-\},\qquad
r\in[7]\setminus\{p\},
\]
where $p^-$ and $q^+$ are read cyclically modulo $5$. For example,
the binomials with $p=1$ are
\[
\begin{aligned}
F_{1,2,r}&=x_2^{(1)}x_1^{(r)}-x_3^{(3)}x_5^{(r)},\\
F_{1,3,r}&=x_3^{(1)}x_1^{(r)}-x_4^{(4)}x_5^{(r)},\\
F_{1,4,r}&=x_4^{(1)}x_1^{(r)}-x_5^{(5)}x_5^{(r)},
\end{aligned}
\qquad r\in\{2,3,4,5,6,7\}.
\]
For each $p$, there are three choices of $q$ and six choices of $r$.
Thus there are $5\cdot3\cdot6=90$ additional binomials.}
\end{Example}

%=%=%=%=%=%=%=%=%=%=%=%=%=%=%=%=%=%=%=%=%=%=%=%=

Set ${\rm in}_{<_{\lex}}(\Bc) := \{{\rm in}_{<_{\lex}}(u) : u \in \Bc\}.$ 

\begin{Lemma}
Every monomial belonging to ${\rm in}_{<_{\lex}}(\Bc)$ is squarefree.   
\end{Lemma}

%A monomial $u=\prod_{i=1}^{2n} x_i^{p_i} \prod_{i=1}^{2n} y_i^{q_i}$ of $S$ is called {\em standard} if $u \not\in {\rm in}_{<_{\lex}}(\Gc)$.  Let $\cM$ denote the set of standard monomials of $S$.

A monomial $u$ of $S$ is called {\em standard} if $u$ is not divisible by any monomial in ${\rm in}_{<_{\lex}}(\Bc)$.  Let $\cM$ denote the set of standard monomials of $S$.

\begin{Lemma}
    If $u \in \cM, v \in \cM$ and $u \neq v$, then $\pi(u) \neq \pi(v)$.
\end{Lemma}

\begin{proof}
Let $\pi(u)=\pi(v)$ and suppose  $\operatorname{supp}(u)\cap\operatorname{supp}(v)=\emptyset$. If only diagonal variables $x_k^{(k)}$ occur in $\operatorname{supp}(u)\cup\operatorname{supp}(v)$, the conclusion follows immediately from $\pi(x_k^{(k)})=a_k$. Otherwise choose the largest variable $x_j^{(i)}$ belonging to $\operatorname{supp}(u)\cup\operatorname{supp}(v)$, where $i\ne j$.  Let $x_j^{(i)} \in \supp(u)$.  Every variable $x_q^{(p)}$ dividing $v$ is strictly smaller than $x_j^{(i)}$. Let $t\ge1$ be the exponent of $x_j^{(i)}$ in $u$.

Given a standard monomial $w$ with $$\pi(w)= \prod_{k=1}^{n} a_k^{\xi_k}\prod_{r=n+1}^{m-1} t_r^{\eta_r}$$ and nonnegative integers $b_1,\ldots,b_n, c_{n+1}, \ldots, c_{m-1}$, 
set
\[
L(w):=\sum_{r=n+1}^{m-1}c_r \eta_r-\sum_{k=1}^{n}b_k \xi_k.
\]
It follows from (\ref{abcdefg}) that
\begin{equation}\label{eq:weight}
L\bigl(x_q^{(p)}\bigr)=b_{p+1}-b_p-b_{q+1},\quad 
L(x_j^{(r)})=c_r - b_{j+1}, \quad
L(x_k^{(k)})=-b_k,
\end{equation}
where $1\le p\le n$, $n<r<m$ and $b_{n+1}=b_1$.  Our job is to choose the coefficients satisfying
\begin{equation}\label{eq:target}
L\bigl(x_j^{(i)}\bigr)=1,\qquad
L(z)=0\ \,\,\,  \text{if} \,\,\, z\in\operatorname{supp}(u)\setminus\{x_j^{(i)}\},\qquad
L(z)\le0\  \, \,\,  \text{if} \,\,\, z<x_j^{(i)}.
\end{equation}

\smallskip

\noindent
$\bullet$\quad Let $i>n$.
Let $c_i=1$ and all other $c_r=0$.  Let $b_k=1$ for $2\le k\le j$
and all other $b_k=0$.  One has $L(x_j^{(i)})=1$.  Since
\[
\supp(u) \setminus\{x_j^{(i)}\} \subset 
\{x_q^{(i)} : 1 \leq q < j\} \cup \{x_q^{(r)} : n < r < i, j \leq q \leq n \} \cup \{ x_{j+1}^{(j+1)}, \ldots, x_{n}^{(n)}, x_1^{(1)}\},
\]
the desired condition (\ref{eq:target}) is satisfied.

\smallskip

\noindent
$\bullet$\quad   Let $i=1$.  Let all $c_r=0$.  Let $b_k=1$ for $2\le k\le j$ and all other $b_k=0$.  One has $L(x_j^{(1)})=1$.  Since
\[
\supp(u) \setminus\{x_j^{(1)}\} \subset \{x_2^{(1)}, \ldots, x_{j-1}^{(1)}\} \cup \{x_1^{(1)}, x_{j+1}^{(j+1)}, \ldots, x_{n}^{(n)}\} \cup \{x_q^{(r)} : n < r < m, j \leq q \leq n \}
\]
the desired condition (\ref{eq:target}) is satisfied.

\smallskip

\noindent
$\bullet$\quad   Let $j=1$ and $1 < i \leq n$.  Let all $c_r=0$.  Let $b_{i+1}=1$ (cyclically) and all other $b_k = 0$.  One has $L(x_1^{(i)})=1$.  Let $i < n$.   Since none of $$x_i^{(1)}, \ldots, x_{i}^{(i-1)},x_{i}^{(n+1)}, \ldots, x_{i}^{(m-1)}, x_{i+1}^{(i+1)}$$ belongs to $\supp(u)$, the desired condition (\ref{eq:target}) is satisfied.  Let $i=n$.  Since none of   
$$
x_n^{(1)}, \ldots, x_{n}^{(n-1)}, 
x_n^{(n+1)}, \ldots, x_{n}^{(m-1)},x_2^{(1)}, \ldots, x_{n-1}^{(1)}, x_1^{(1)}
$$
belongs to $\supp(u)$, the desired condition (\ref{eq:target}) is satisfied.

\smallskip

\noindent
$\bullet$\quad   Let $1 < i \leq n, j > 1$. Given a standard monomial $w$, we introduce
\begin{eqnarray*}
R(w)&:=&\{p:x_q^{(p)}\in \supp(w) \text{ for some }q\ne p\},\\
C(w)&:=&\{q:x_q^{(p)}\in \supp(w) \text{ for some }p\ne q\}.
\end{eqnarray*}
One has $$R(w) \cap C(w) \subset \{1\}.$$  Let 
\[
\varepsilon=
\begin{cases}
1,&n\in C(u),\\
0,&n\notin C(u)
\end{cases}
\]
and set $c_r = \varepsilon$ for $n<r<m$.  Since $j\in C(u)$, it follows that $j\notin R(u)$.  Define the coefficients recursively as follows:
\begin{eqnarray}
b_1&=&1,\label{eq:b1xxx}\\
b_{p+1}&=&\begin{cases}
b_p+\varepsilon, \, \, \, \, \, \, \, \, \, \, \, \, \, \, \, p\in R(u),\\
\varepsilon+\mathbf1_{\{p<j\}}, \, \, \, \, \, p\notin R(u),
\end{cases}\qquad 1\le p<i,\label{eq:bpxxx}\\
b_{i+1}&=&b_i+1+\varepsilon,\label{eq:aaabi}\\
b_k&=&\varepsilon+\mathbf1_{\{k\le j\}},\qquad i+2\le k\le n.\label{eq:bkxxx}
\end{eqnarray}
Here $\mathbf1$ is $1$ if its indicated condition holds and $0$ otherwise. When $i=n$, one has $\varepsilon = 0$ and $j \not\in R(u)$ makes $b_{j+1}=\cdots=b_n=0$, so (\ref{eq:aaabi}) gives $b_{n+1}=1$.
% 
% For clarity, the elementary consequences of the construction that we need are
We claim that
\begin{equation}\label{eq:consequencesxxx}
\begin{gathered}
b_k\ge\varepsilon\quad(1\le k\le n),\qquad b_k\ge1+\varepsilon\quad(2\le k\le j),\\
b_{p+1}-b_p\le\varepsilon\quad(p<i),\qquad b_{j+1}=\varepsilon.
\end{gathered}
\end{equation}
In fact, the first and second inequalities follow along the recursion.  For the third inequality, equality holds at an active row $p\in R(u)$; at an inactive row the next coefficient is $\varepsilon+\mathbf1_{\{p<j\}}$, which is at most $b_p+\varepsilon$.  For the last equality, if $j<i$ use $j\notin R(u)$ in (\ref{eq:bpxxx}); if $i<j<n$ use (\ref{eq:bkxxx}); if $j=n$, then $\varepsilon=1=b_1$. % This explicitly covers both $i<j$ and $j<i$.

  Now, equations \eqref{eq:weight}, (\ref{eq:aaabi}) and (\ref{eq:consequencesxxx}) show that $L(x_j^{(i)})=1$.  We show that if $x_q^{(p)} < x_j^{(i)}$, then $L(x_q^{(p)}) \leq 0$.  Clearly, $L(x_k^{(k)}) = -b_k \leq 0$.  
 %For $n < r < m $,
 % \[
 % L(x_j^{(r)}) = \varepsilon - b_{j+1} \leq 0.
 % \]
 For $n<r<m$ and $1\le q\le n$,
\[
L(x_q^{(r)})=\varepsilon-b_{q+1}\le0.
\]
For $p<i$,
\[
L\bigl(x_q^{(p)}\bigr)
= b_{p+1}-b_{p}-b_{q+1} \le\varepsilon-b_{q+1}\le0;
\]
and in row $i$, for $q<j$ with $q\ne i$,
\[
L\bigl(x_q^{(i)}\bigr)
= b_{i+1}-b_{i}-b_{q+1}
= 1+\varepsilon-b_{q+1}\le0.
\]
It remains to check that every other factor $x_q^{(p)}\in  \supp(u)$ satisfies $L(x_q^{(p)})=0$.
\begin{enumerate}
\item[(i)]
Let $x_q^{(p)}\in \supp(u)$ with $p \neq q$ and either $p<i$ or $p > n$.  One has $p\in R(u)$. If $p<i$  then $b_{p+1}-b_p=\varepsilon$ and if $p>n$ then $c_p=\epsilon$. One has $q\ge j > 1, q\ne i$ and $q\notin R(u)$.  If $q<i$, then (\ref{eq:bpxxx}) says that $b_{q+1} = \varepsilon$; if $i<q<n$, then (\ref{eq:bkxxx}) says that $b_{q+1} = \varepsilon$; and if $q=n$, then $b_{q+1} = b_1 = \varepsilon$.  So,  $L(x_q^{(p)})=b_{p+1}-b_p-b_{q+1}=0$ for $p <i$ and $L(x_q^{(p)})=\varepsilon - b_{q+1}=0$ for $p > n$.

\item[(ii)] Let $x_q^{(i)}\in \supp(u)$ with $q \neq i$ and $q<j$.  One has $b_{q+1}=1+\varepsilon$.  In fact, for $q=1$, this follows from $b_1=1$ and (\ref{eq:bpxxx}), whether or not $1\in R(u)$. If $q>1$, then $q\notin R(u)$, so the assertion follows from (\ref{eq:bpxxx}) or (\ref{eq:bkxxx}). 
\item[(iii)] Let $x_k^{(k)}\in \supp(u)$.
One has $k>j$ and $k\ne i+1$ with the latter index read cyclically.  % For $2\le k\le j$, use the minor on row $i$, the last row, and columns $k-1,j$; its leading monomial is $d_kx_j^{(i)}$. For $k=1$, use $d_1x_j^{(i)}-x_j^{(1)}x_1^{(i)}$. For $k=i+1$, use $d_id_{j+1}-x_j^{(i)}d_{i+1}$. 
Furthermore, one has $n\notin C(u)$.  So, $\varepsilon=0$.  %Since $j \not\in R(u)$, $b_{j+1} = \cdots = b_i = 0$ follows from (\ref{eq:bpxxx}).  In addition, $b_{i+2} = \cdots = b_n = 0$ follows from (\ref{eq:bkxxx}).
If $j<i$, then $j\notin R(u)$ and (\ref{eq:bpxxx}) give
$b_{j+1}=\cdots=b_i=0$.
If $k\ge i+2$, then $k>j$ and (\ref{eq:bkxxx}) give $b_k=0$.
Hence $L(x_k^{(k)})=-b_k=0$.
\end{enumerate}
This completes the verification of \eqref{eq:target}.

Finally, additivity gives $L(u) = t > 0$ and $L(v) \leq 0$.  This contradicts $\pi(u)=\pi(v)$.
\hspace{0.9cm}
\end{proof}

A standard fact on Gr\"obner bases
\cite[Proposition~2.2.5]{HHgtm260},
together with the normality criterion
\cite[Proposition~13.15]{Stu}
and Hochster's theorem
\cite[Theorem~1]{Hochster1972},
guarantees that

\begin{Corollary}
The set of quadratic binomials belonging to $\Bc$
is a Gr\"obner basis of $\Ker(\pi)$ with respect
to $<_{\lex}$, and its initial monomials are
squarefree. In particular, $I(X)=\Ker(\pi)$
and $S/I(X)$ is a normal and Cohen--Macaulay domain.
\end{Corollary}

This completes the proof of Theorem~\ref{thm:prime-classification}.

%%%%%%%%%%%%%%%%%%%%%%%%%%%%%%%%%%%%%%%%%%%%%%%%%%%%%%%%%%

\section{Regularity and Graded Betti Numbers of $I(X_{m,n})$}\label{sec:3}

We determine the graded Betti numbers and regularity of the prime
ideals in Theorem~\ref{thm:prime-classification}. For $m\ge4$, the quotient has the
same graded Betti numbers as the corresponding generic determinantal
ring. For $m=3$, the scroll presentation gives the Betti numbers by the
Eagon--Northcott resolution. We then prove that the squarefree
quadratic initial ideal constructed in
Section~\ref{subsec:m=3} has a linear resolution and that the
degeneration to this ideal preserves the graded Betti numbers.

%%%%%%%%%%%%%%%%%%%%%%%%%%%%%%%%%%%%%%%%%%%%%%%%%%%%%%%%%%%%%

\subsection{The case $m\geq 4$}

We begin with the following consequence of the results of Section~2.

\begin{Lemma}\label{lem:dimension}
Assume that $m\geq 4$ and that $\mathcal G_X$ consists of a single
row-distinct directed cycle of length $n$. Then
\[
\dim S/I(X_{m,n})=m-1.
\]
\end{Lemma}

\begin{proof}
By Corollary~2.9, we have
$S/I(X_{m,n})\cong\operatorname{im}\pi$, where
$\pi$ is the toric parametrization defined in
Section~2.2. By (\ref{abcdefg}), the fraction field of
$\operatorname{im}\pi$ has transcendence degree $m-1$
over $K$. Consequently,
$\dim S/I(X_{m,n})=m-1$.
\end{proof}

%=%=%=%=%=%=%=%=%=%=%=%=%=%=%=%=%=%=%=%=%=%=%=%=%=%=%=%=%=%=
We now determine the graded Betti numbers of
$S/I(X_{m,n})$ for $m\geq 4$. As a consequence,
we obtain the regularity of $I(X_{m,n})$ and
complete the proof of Theorem~\ref{thm:PCML-classification}(ii).

\begin{Proposition}\label{prop:regularity}
Assume that $m\geq 4$ and that $\mathcal G_X$ consists of a
single row-distinct directed cycle of length $n$.
Let $Z$ be a generic $m\times n$ matrix of variables
over $K$, and let $T$ be the polynomial ring generated
by its entries. Then
\[
\beta_{i,j}^{S}(S/I(X_{m,n}))
=
\beta_{i,j}^{T}(T/I_2(Z))
\]
for all $i,j$. In particular,
\[
\operatorname{reg} I(X_{m,n})=n
\quad\text{and}\quad
\operatorname{pd}_S(S/I(X_{m,n}))=(m-1)(n-1).
\]
Consequently, $I(X_{m,n})$ has a linear resolution if and only if $n=2$.
\end{Proposition}

\begin{proof}
Set $I=I(X_{m,n})$ and $R=S/I$. Since $\mathcal G_X$ consists of a single row-distinct directed
cycle of length $n$, we have $n\leq m-1$.
After simultaneously relabeling the columns and the first $m-1$
rows, we may assume that $X_{m,n}$ is the matrix considered in
Section~\ref{sec:m=4}, so that
\begin{equation}\label{eq:14}
y_i=x_{i+1}^{(i+1)}
\quad\text{for }1\leq i<n,
\qquad
y_n=x_1^{(1)}.
\end{equation}

Let $T=K\big[z_j^{(i)}:1\leq i\leq m,\ 1\leq j\leq n\big]$
and let $Z=\big(z_j^{(i)}\big)$ be the generic $m\times n$ matrix.
Set $A=T/I_2(Z)$. It is well known that the generic determinantal
ring $A$ is a Cohen--Macaulay domain with $\dim A=m+n-1$ and $\operatorname{reg} A=\min\{m,n\}-1$; see, for example,
\cite{BrunsVetter} or \cite{BrunsConcaRaicuVarbaro}.
Since $n\leq m-1$, we have $\operatorname{reg} A=n-1$.
Consider in $A$ the linear forms
\[
\ell_i=z_i^{(m)}-z_{i+1}^{(i+1)}
\quad\text{for }1\leq i<n,
\quad\text{and}\quad
\ell_n=z_n^{(m)}-z_1^{(1)}.
\]
In view of~\eqref{eq:14}, imposing these $n$ linear relations on
the generic matrix $Z$ gives precisely the matrix $X_{m,n}$. Hence $R\cong A/(\ell_1,\ldots,\ell_n).$ It follows from Lemma~\ref{lem:dimension} that $\dim R=m-1$.
Thus, we have
$\dim A/(\ell_1,\ldots,\ell_n)=m-1=\dim A-n$.
Since $A$ is Cohen--Macaulay, it follows that
$\ell_1,\ldots,\ell_n$ is an $A$-regular sequence; see, for example,
\cite[Section~2.1]{BrunsHerzog}.
Moreover, these linear forms are linearly independent in $T$,
so they also form a regular sequence on $T$, and hence
$S\cong T/(\ell_1,\ldots,\ell_n)$.
Therefore, tensoring a minimal graded free resolution of $A$
over $T$ with $S$ gives a minimal graded free resolution of $R$
over $S$. Consequently,
$\beta_{i,j}^{S}(R)=\beta_{i,j}^{T}(A)$ for all $i,j$.
In particular, $R$ and $A$ have the same regularity and
projective dimension. Hence
\[
\operatorname{reg} I(X_{m,n})=n
\quad\text{and}\quad
\operatorname{pd}_S(R)=(m-1)(n-1).
\]
Since $I(X_{m,n})$ is generated by quadratic forms,
it follows that $I(X_{m,n})$ has a linear resolution
if and only if $n=2$.
\end{proof}
%=%=%=%=%=%=%=%=%=%=%=%=%=%=%=%=%=%=%=%=%=%=%=%=%=%=%=%=%=%=%=%=%=%=%=%=%=%=%=%=%=%=%=%=%=%=%=%=%=%=%=%=

\subsection{The case $m=3$} \label{subec:3.2}

Assume that $I(X_{3,n})$ is prime and write $n=2q$.
Section~\ref{subsec:m=3} identifies $I(X_{3,n})$ with the
ideal of $2$-minors of a $2\times3q$ matrix consisting of $q$ scroll
blocks of size three. We first record the graded Betti numbers
arising from this presentation.

\begin{Corollary}\label{cor:nd-scroll-betti}
Assume that $m=3$ and that $\mathcal G_X$ is a disjoint
union of row-distinct directed cycles of length $2$.
Write $n=2q$. Then
\[
\beta_{i,i+1}^{S}(S/I(X_{3,n}))
=i\binom{3q}{i+1},
\qquad 1\leq i\leq 3q-1,
\]
with $\beta_{0,0}^{S}(S/I(X_{3,n}))=1$,
and all other graded Betti numbers vanish.
\end{Corollary}

\begin{proof}
By Section~\ref{subsec:m=3}, $I(X_{3,n})=I_2(H)$,
where $H$ is a $2\times3q$ scroll matrix. By
\cite[Corollary~A2.62]{EisenbudSyzygies},
$I_2(H)$ has height $3q-1$, and the Eagon--Northcott complex
gives a minimal graded free resolution of $S/I_2(H)$.
Its module in homological degree $i$ is
$S(-i-1)^{\,i\binom{3q}{i+1}}$ for $1\le i\le3q-1$.
The assertions follow.
\end{proof}

Together with Proposition~\ref{prop:regularity}, this proves
Theorem~\ref{thm:PCML-classification}.

%=%=%=%=%=%=%=%=%=%=%=%=%=%=%=%=%=%=%=%=%=%=%=%=%=%=%=%=%=%=%=%=

We now consider the initial ideal arising from the lexicographic
order used in Section~\ref{subsec:m=3}. After relabeling
the columns and variables, write
\[
X_{3,n}=
\begin{pmatrix}
x_1&\cdots&x_q&x_{q+1}&\cdots&x_{2q}\\
y_1&\cdots&y_q&y_{q+1}&\cdots&y_{2q}\\
y_{q+1}&\cdots&y_{2q}&x_1&\cdots&x_q
\end{pmatrix}.
\]
Let $<_{\mathrm{lex}}$ be the lexicographic order on
$S=K[x_1,\ldots,x_{2q},y_1,\ldots,y_{2q}]$ induced by
\[
x_{q+1}>\cdots>x_{2q}>y_1>\cdots>y_q>
x_1>\cdots>x_q>y_{q+1}>\cdots>y_{2q},
\]
and set
\[
J=\operatorname{in}_{<_{\mathrm{lex}}}(I(X_{3,n})).
\]
By Section~\ref{subsec:m=3}, $J$ is generated by the
squarefree initial monomials of the $2$-minors of $X_{3,n}$.
We prove that $J$ has a linear resolution by exhibiting a perfect
elimination ordering for the complement of its associated graph.
We first recall the relevant graph-theoretic definitions and
Fr\"oberg's theorem.

Let $G$ be a finite simple graph, with vertex set $V(G)$ and edge set $E(G)$, and let $\overline{G}$ be its complement. A subset of $V(G)$ is called a \emph{clique} if its elements are pairwise adjacent, and $G$ is called \emph{chordal} if it has no induced cycle of length at least $4$. Given an ordering $(v_1,\ldots,v_N)$ of $V(G)$, the \emph{later neighborhood} $N_G^+(v_i)$ of $v_i$ is the set of the neighbors of $v_i$ in $\{v_{i+1},\ldots,v_N\}$. The ordering is called a \emph{perfect elimination ordering} if $N_G^+(v_i)$ is a clique for every $i=1,\ldots,N$. It is well known that a graph is chordal if and only if it admits a perfect elimination ordering \cite{FG,Rose}.

If the vertices of a finite simple graph $G$ are identified with the
variables of a polynomial ring, its edge ideal is
\[
I(G)=(uv:\{u,v\}\in E(G)).
\]
We recall Fr\"oberg's theorem \cite{F}, which states that
$I(G)$ has a linear resolution if and only if the complement
$\overline G$ is chordal.

\begin{Proposition}\label{prop:m=3linres}
With the notation above, the initial ideal $J=\operatorname{in}_{<_{\mathrm{lex}}}(I(X_{3,n}))$ has a linear resolution.
\end{Proposition}

\begin{proof}
Since $J$ is generated by squarefree monomials of degree $2$, it is
	the edge ideal of a finite simple graph $G$. Comparing the two terms of each $2$-minor of
$X_{3,n}$ with respect to $<_{\lex}$, one obtains
\[		I(G)=J
		=\bigl(
		x_{q+i}y_j,\,
		x_{q+i}y_{q+j},\,
		y_ix_j
		:1\leq i,j\leq q
		\bigr)+
		\bigl(
		x_{q+i}x_j,\,
		y_i y_{q+j},\,
		y_{q+j}x_i
		:1\leq i<j\leq q
		\bigr).
\]
	
	We prove that the complement graph $\overline{G}$ is chordal. 
    We begin by describing its edges; the reader can look at Figure \ref{Figure: complement graph}. Let us set $V_1=\{x_{q+1},\ldots,x_{2q}\}$, $V_2=\{y_1,\ldots,y_q\}$, $V_3=\{x_1,\ldots,x_q\}$ and $V_4=\{y_{q+1},\ldots,y_{2q}\}$.
    No generator of $J$ is the product of two variables of the same
$V_i$, so $V_i$ is a clique of $\overline{G}$, for all $i\in [4]$. Since
$x_{q+i}y_j$, $x_{q+i}y_{q+j}$ and $y_ix_j$ belong to $J$ for all
$1\leq i,j\leq q$, there are no edges between $V_1$ and $V_2$,
between $V_1$ and $V_4$, or between $V_2$ and $V_3$. Moreover, observe that the monomials $x_{q+i}x_j$, $y_iy_{q+j}$ and
$x_iy_{q+j}$ belong to $J$ exactly when $i<j$; hence, for
$1\leq i,j\leq q$, each of
$\{x_{q+i},x_j\}$, $\{y_i,y_{q+j}\}$ and $\{x_i,y_{q+j}\}$ is an edge of
$\overline{G}$ if and only if $i\geq j$. The preceding observations describe all the edges of
	$\overline{G}$.

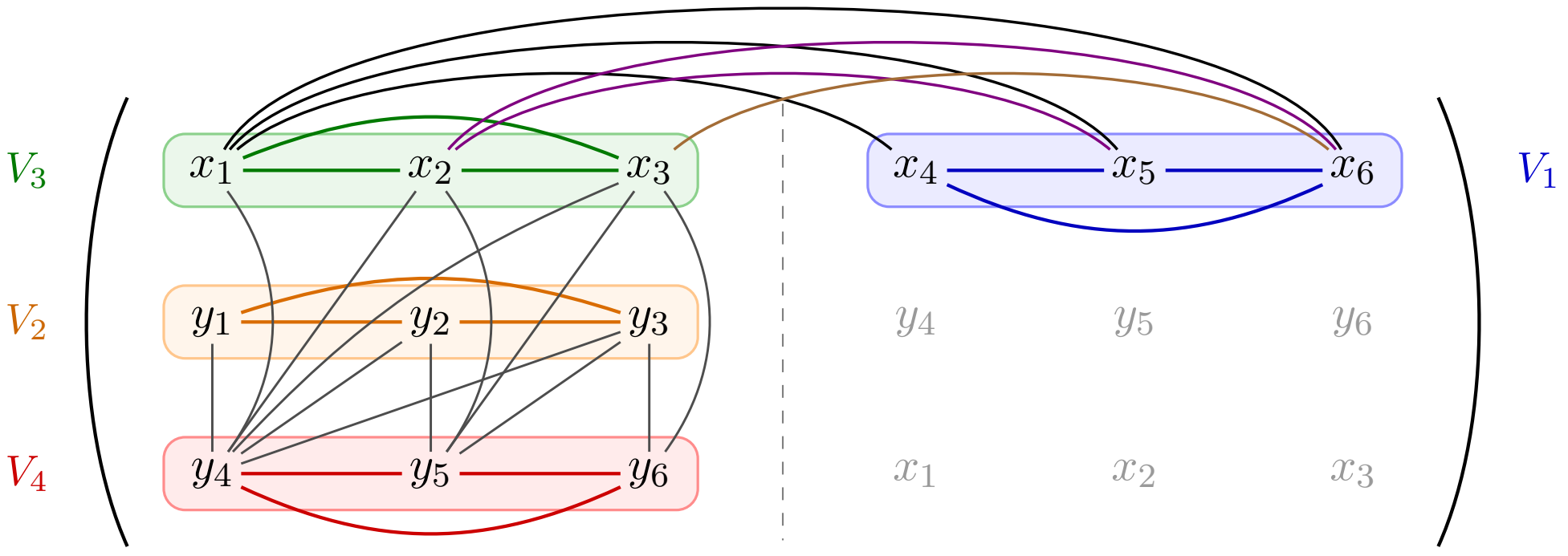
\begin{figure}[h]
\centering
\begin{tikzpicture}[x=1cm, y=1cm,
  entry/.style={inner sep=1.5pt, font=\small},
  copy/.style={inner sep=1.5pt, font=\small, text=black!40},
  clique/.style={rounded corners=5pt, fill=#1!8, draw=#1!45, line width=0.6pt},
  cliquelab/.style={font=\footnotesize, text=#1!80!black},
  inedge/.style={line width=0.8pt, draw=#1},
  xedge/.style={line width=0.55pt, draw=black!70},
  xarc/.style={line width=0.65pt, draw=#1}
]
% Colors of the four cliques
\colorlet{Vone}{blue}
\colorlet{Vtwo}{orange}
\colorlet{Vthree}{green!60!black}
\colorlet{Vfour}{red}
 
% Cliques: V_3, V_2, V_4 in the left block, V_1 in the first row on the right
\draw[clique=Vthree] (-0.4, 0.3)  rectangle (4.0,-0.3);   % V_3
\draw[clique=Vtwo]   (-0.4,-0.95) rectangle (4.0,-1.55);  % V_2
\draw[clique=Vfour]  (-0.4,-2.2)  rectangle (4.0,-2.8);   % V_4
\draw[clique=Vone]   ( 5.4, 0.3)  rectangle (9.8,-0.3);   % V_1
 
% Separator between the two blocks
\draw[dashed, gray] (4.7,0.55) -- (4.7,-3.05);
 
% Parentheses
\draw[line width=0.8pt]
  (-0.7,0.6) .. controls (-1.15,-0.4) and (-1.15,-2.1) .. (-0.7,-3.1);
\draw[line width=0.8pt]
  (10.1,0.6) .. controls (10.55,-0.4) and (10.55,-2.1) .. (10.1,-3.1);
 
% Matrix entries: row 1 = a1..a6, row 2 = b1..b6, row 3 = c1..c6
\foreach \idx/\col/\lab in {1/0/x_1, 2/1.8/x_2, 3/3.6/x_3, 4/5.8/x_4, 5/7.6/x_5, 6/9.4/x_6}
  \node[entry] (a\idx) at (\col,0) {$\lab$};
\foreach \idx/\col/\lab in {1/0/y_1, 2/1.8/y_2, 3/3.6/y_3}
  \node[entry] (b\idx) at (\col,-1.25) {$\lab$};
\foreach \idx/\col/\lab in {1/0/y_4, 2/1.8/y_5, 3/3.6/y_6}
  \node[entry] (c\idx) at (\col,-2.5) {$\lab$};
% repeated variables (gray copies, not new vertices)
\foreach \col/\lab in {5.8/y_4, 7.6/y_5, 9.4/y_6}
  \node[copy] at (\col,-1.25) {$\lab$};
\foreach \col/\lab in {5.8/x_1, 7.6/x_2, 9.4/x_3}
  \node[copy] at (\col,-2.5) {$\lab$};
 
% Edges inside the cliques (each clique is a K_3, drawn once;
% the arc of V_1 goes below the first row to leave room for E_{13})
\draw[inedge=Vthree!80!black] (a1) -- (a2) -- (a3);
\draw[inedge=Vthree!80!black] (a1) to[bend left=22] (a3);
\draw[inedge=Vtwo!85!black]   (b1) -- (b2) -- (b3);
\draw[inedge=Vtwo!85!black]   (b1) to[bend left=18] (b3);
\draw[inedge=Vfour!80!black]  (c1) -- (c2) -- (c3);
\draw[inedge=Vfour!80!black]  (c1) to[bend right=25] (c3);
\draw[inedge=Vone!75!black]   (a4) -- (a5) -- (a6);
\draw[inedge=Vone!75!black]   (a4) to[bend right=25] (a6);
 
% E_{34}: rows 1 and 3 of the left block, region of the T_{ij}
% (the bent edges avoid the entries of the second row)
\draw[xedge] (a1) to[bend left=36]  (c1);  % x_1 y_4
\draw[xedge] (a2) --                (c1);  % x_2 y_4
\draw[xedge] (a2) to[bend left=36]  (c2);  % x_2 y_5
\draw[xedge] (a3) to[bend right=12] (c1);  % x_3 y_4
\draw[xedge] (a3) --                (c2);  % x_3 y_5
\draw[xedge] (a3) to[bend left=36]  (c3);  % x_3 y_6
 
% E_{24}: rows 2 and 3 of the left block, region of the Q_{ij}
\draw[xedge] (b1) -- (c1);  % y_1 y_4
\draw[xedge] (b2) -- (c1);  % y_2 y_4
\draw[xedge] (b2) -- (c2);  % y_2 y_5
\draw[xedge] (b3) -- (c1);  % y_3 y_4
\draw[xedge] (b3) -- (c2);  % y_3 y_5
\draw[xedge] (b3) -- (c3);  % y_3 y_6
 
% E_{13}: arcs over the first row, x_j -- x_{3+i} for j <= i,
% colored by their endpoint in V_3
\draw[xarc=black]         (a1) to[bend left=40, looseness=0.62] (a4);  % x_1 x_4
\draw[xarc=black]         (a1) to[bend left=50, looseness=0.54] (a5);  % x_1 x_5
\draw[xarc=black]         (a1) to[bend left=60, looseness=0.50] (a6);  % x_1 x_6
\draw[xarc=violet]        (a2) to[bend left=40, looseness=0.62] (a5);  % x_2 x_5
\draw[xarc=violet]        (a2) to[bend left=50, looseness=0.54] (a6);  % x_2 x_6
\draw[xarc=brown!85!black](a3) to[bend left=40, looseness=0.62] (a6);  % x_3 x_6
 
% Names of the cliques
\node[cliquelab=Vthree, anchor=east] at (-1.2, 0)    {$V_3$};
\node[cliquelab=Vtwo,   anchor=east] at (-1.2,-1.25) {$V_2$};
\node[cliquelab=Vfour,  anchor=east] at (-1.2,-2.5)  {$V_4$};
\node[cliquelab=Vone,   anchor=west] at (10.6, 0)    {$V_1$};
\end{tikzpicture}
\caption{The graph $\overline{G}$ for $q=3$.}
\label{Figure: complement graph}
\end{figure}

	Now let $\mathcal O$ be the following ordering of vertices of $\overline{G}$: 
    $$\mathcal{O}=(x_{q+1},\ldots,x_{2q},\;
		y_1,\ldots,y_q,\;
		x_1,\ldots,x_q,\;
		y_{q+1},\ldots,y_{2q}).$$ 
    We claim that $\mathcal{O}$ is a perfect elimination ordering of $\overline{G}$. 
	
	Fix $1\leq i\leq q$. From the description of the edges of $\overline{G}$ it is easy to see that:
    \begin{itemize}
        \item $
	N_{\overline{G}}^+(x_{q+i})
	=
	\{x_{q+k}:i<k\leq q\}
	\cup
	\{x_j:1\leq j\leq i\};$
    \item $N_{\overline{G}}^+(y_i)
	=
	\{y_k:i<k\leq q\}
	\cup
	\{y_{q+j}:1\leq j\leq i\};$
    \item $N_{\overline{G}}^+(x_i)
	=
	\{x_k:i<k\leq q\}
	\cup
	\{y_{q+j}:1\leq j\leq i\}$;
    \item $N_{\overline{G}}^+(y_{q+i})
	=
	\{y_{q+k}:i<k\leq q\}.$
    \end{itemize}
    The last
one is contained in $V_4$, hence it is a clique. Each of the other
three is the union of a subset of one clique and a subset of another
one (of $V_1$ and $V_3$, of $V_2$ and $V_4$, and of $V_3$ and $V_4$,
respectively), so it suffices to check that every vertex of the first
subset is adjacent to every vertex of the second. Such a pair has the
form $\{x_{q+k},x_j\}$, $\{y_k,y_{q+j}\}$ or $\{x_k,y_{q+j}\}$ with
$k>i\geq j$, and it is an edge since $k\geq j$.
Hence $\mathcal O$ is a perfect elimination ordering of
$\overline{G}$, as we claimed. 

It follows that $\overline{G}$ is chordal and, in conclusion, that $J=\operatorname{in}_{<_{\mathrm{lex}}}(I(X_{3,n}))$ has a linear resolution. This concludes the proof.
\end{proof}

%=%=%=%=%=%=%=%=%=%=%=%=%=%=%=%=%=%=%=%=%=%=%=%=%=%=%=%=%=

    \begin{Corollary}\label{cor:m3-linear}
With the notation above, one has
\[
\beta^S_{i,j}(S/I(X_{3,n}))
=
\beta^S_{i,j}(S/J)
\]
for all $i,j$. In particular,
\[
\beta^S_{i,i+1}(S/J)
=
i\binom{3q}{i+1},
\qquad 1\le i\le3q-1,
\]
with $\beta^S_{0,0}(S/J)=1$, and all other graded Betti numbers vanish.
\end{Corollary}

\begin{proof}
By Corollary~\ref{cor:nd-scroll-betti}, the ideal $I(X_{3,n})$ has a
$2$-linear resolution. By Proposition~\ref{prop:m=3linres},
the same holds for $J$. For either $L=I(X_{3,n})$ or $L=J$, the
Hilbert series is therefore
\[
\operatorname{Hilb}_{S/L}(t)
=
\frac{1+\sum_{i\ge1}(-1)^i
\beta^S_{i,i+1}(S/L)t^{i+1}}{(1-t)^{4q}}.
\]
Since an ideal and its initial ideal have the same Hilbert series,
the numerators agree. Comparing the coefficient of $t^{i+1}$ gives
the equality of the graded Betti numbers. The displayed formula
then follows from Corollary~\ref{cor:nd-scroll-betti}.
\end{proof}

%=%=%=%=%=%=%=%=%=%=%=%=%=%=%=%=%=%=%=%=%=%=%=%=%=%=%=%=%=
%=%=%=%=%=%=%=%=%=%=%=%=%=%=%=%=%=%=%=%=%=%=%=%=%=%=%=%=%=

\begin{Remark}\label{rem:pcml-scrolls}
{\em
The PCML matrices in Theorem~\ref{thm:PCML-classification} have two-row presentations whose scroll types can be read from the path components described
in Appendix~\ref{app:two-row}.

If $m=3$, write $n=2q$. The matrix $H$ constructed in
Section~\ref{subsec:m=3} consists of $q$ scroll blocks of
size three. Its directed graph has $q$ path components, each on
three vertices. Hence $I(X_{3,n})$ defines the scroll
$S(3,\ldots,3)$, with $q$ entries.

If $n=2$, write
\[
y_1=x_2^{(i)},\qquad y_2=x_1^{(k)},
\qquad i\ne k.
\]
After ordering the columns of the transpose as the $i$-th row,
the last row, the $k$-th row, and the remaining rows of $X_{m,2}$,
we obtain
\[
\begin{pmatrix}
x_1^{(i)}&x_2^{(i)}&x_1^{(k)}&
x_1^{(r_1)}&\cdots&x_1^{(r_{m-3})}\\
x_2^{(i)}&x_1^{(k)}&x_2^{(k)}&
x_2^{(r_1)}&\cdots&x_2^{(r_{m-3})}
\end{pmatrix},
\]
where $\{r_1,\ldots,r_{m-3}\}=[m-1]\setminus\{i,k\}$.
The first three columns form a scroll block of size three,
and each remaining column forms a scroll block of size one.
Thus the directed graph consists of one path on three vertices
and $m-3$ isolated vertices. Transposition preserves the ideal of
$2$-minors, so $I(X_{m,2})$ defines the scroll
$S(3,1,\ldots,1)$, with $m-3$ entries equal to one.
}
\end{Remark}

%=%=%=%=%=%=%=%=%=%=%=%=%=%=%=%=%=%=%=%=%=%=%=%=%=%=%=%=%=%=%=%=%=%=%=

The normality of every prime quotient in our classification was
established in Section~\ref{sec:2}. For the PCML cases, the
following remark relates normality to the classical theory of
varieties of minimal degree.

\begin{Remark}
{\em 
Let $I\subseteq S=K[x_1,\ldots,x_N]$ be a homogeneous
prime ideal generated by quadratic forms such that
$S/I$ is Cohen--Macaulay and $I$ has a linear
resolution. If $K$ is algebraically closed, then
$S/I$ is normal. Indeed, $I$ defines a variety of
minimal degree, and the classical classification
of such varieties implies normality; see
\cite{EG84,EH87}.

Over an arbitrary field, however, these conditions
do not necessarily imply normality. For example,
$I=(x^2+y^2)\subseteq\mathbb{R}[x,y]$ is prime,
$S/I$ is Cohen--Macaulay, and $I$ has a linear
resolution, but $S/I$ is not normal.

On the other hand, if $I$ is a toric ideal, then
$S/I$ is an affine semigroup ring. Since a toric
ideal remains prime after extending the field
to its algebraic closure, and the normality of
an affine semigroup ring is independent of the
base field \cite[Proposition~1]{Hochster1972},
it follows that $S/I$ is normal over any field $K$.
}
\end{Remark}

%=%=%=%=%=%=%=%=%=%=%=%=%=%=%=%=%=%=%=%=%=%=%=%=%=%=%=%=%=%=%=%=

\appendix
\section{Two-row matrices and rational normal scrolls}
\label{app:two-row}
In this appendix, we characterize the two-row matrices whose
ideals of $2$-minors are prime in terms of directed paths and
identify these ideals with the defining ideals of rational
normal scrolls.

Let $X$ be a $2\times n$ matrix of not necessarily distinct
variables over a field $K$, and let $S$ be the polynomial ring
generated by the distinct variables occurring in $X$. We assume $I_2(X)\ne0$ and delete repeated columns.
Since primality forces the entries in each row and in each
column to be pairwise distinct, we restrict attention to matrices
\begin{equation}\label{eq:app-two-row-matrix}
X=
\begin{pmatrix}
x_1&x_2&\cdots&x_n\\
z_1&z_2&\cdots&z_n
\end{pmatrix},
\end{equation}
where the $x_j$ are pairwise distinct, the $z_j$ are pairwise
distinct, and $z_j\ne x_j$ for every $j\in[n]$. We associate with $X$ the directed graph $\mathcal G_X$ on $[n]$ in which $j\longrightarrow i$ if and only if $z_j=x_i.$ Thus the arrow points from the column containing the repeated entry in
the second row to the column containing that variable in the first row,
as in the definition of the graph of a nearly determinantal matrix.
Every vertex has indegree and outdegree at most one, and there are no
loops. Hence every connected component of $\mathcal G_X$ is a directed
path or a directed cycle. Here an isolated vertex is regarded as a
directed path, and a directed cycle has length at least two.

%=%=%=%=%=%=%=%=%=%=%=%=%=%=%=%=%=%=%=%=%=%=%=%=%=%=%=%=%=%=%=%=

\begin{Lemma}\label{lem:app-cycle}
Let $r\ge2$, let $R=K[u_1,\ldots,u_r]$, and set
\[
U=
\begin{pmatrix}
u_1&u_2&\cdots&u_r\\
u_2&u_3&\cdots&u_1
\end{pmatrix}.
\]
Then $I_2(U)$ is not prime.
\end{Lemma}

\begin{proof}
Put $u_{r+1}=u_1$. Summing the minors determined by the first column and
each of the remaining columns gives $\sum_{j=2}^{r}(u_1u_{j+1}-u_2u_j)
=(u_1-u_2)(u_1+\cdots+u_r)\in I_2(U).$ Both factors are nonzero linear forms, whereas $I_2(U)$ contains no
nonzero linear form. Hence $I_2(U)$ is not prime.
\end{proof}

For positive integers $r_1,\ldots,r_s$, consider the matrices
\begin{equation}\label{eq:app-scroll-blocks}
H=(H_1\mid\cdots\mid H_s),
\qquad
H_j=
\begin{pmatrix}
w_{j,0}&w_{j,1}&\cdots&w_{j,r_j-1}\\
w_{j,1}&w_{j,2}&\cdots&w_{j,r_j}
\end{pmatrix},
\end{equation}
where the variables $w_{j,k}$, with $1\le j\le s$ and
$0\le k\le r_j$, are pairwise distinct. Put $N=r_1+\cdots+r_s$ and
$A=K[w_{j,k}:1\le j\le s,\ 0\le k\le r_j]$, endowed with
the standard grading. The projective variety $S(r_1,\ldots,r_s)
=
\operatorname{Proj}(A/I_2(H))
\subseteq \mathbb P_K^{N+s-1}$ is called the \emph{rational normal scroll of type
$(r_1,\ldots,r_s)$}; see \cite[Appendix~A2H, p.~224]{EisenbudSyzygies}.
Its homogeneous coordinate ring is $A/I_2(H)$.
Each $H_j$ is called a \emph{scroll block of size $r_j$}.
Thus the type records the numbers of columns in the blocks
and is considered up to permutation.
When $s=1$, this construction gives the rational normal
curve of degree $r_1$ in $\mathbb P_K^{r_1}$.

\begin{Theorem}\label{thm:app-two-row-classification}
Let $X$ be as in~\eqref{eq:app-two-row-matrix}, with $n\ge2$. Then the
following conditions are equivalent:
\begin{enumerate}
\item[$\mathrm{(i)}$] $X$ is a PCML matrix;
\item[$\mathrm{(ii)}$] $I_2(X)$ is prime;
\item[$\mathrm{(iii)}$] $\mathcal G_X$ contains no directed cycle;
\item[$\mathrm{(iv)}$] after a permutation of the columns and a relabeling
of the variables, $X$ is the matrix $H$ in~\eqref{eq:app-scroll-blocks}
for some positive integers $r_1,\ldots,r_s$ satisfying
$r_1+\cdots+r_s=n$.
\end{enumerate}
When these conditions hold, the integers $r_1,\ldots,r_s$ are the numbers
of vertices in the connected components of $\mathcal G_X$. In particular,
$I_2(X)$ is identified with the defining ideal of the rational normal
scroll of type $(r_1,\ldots,r_s)$.
\end{Theorem}

\begin{proof}
The implication $\mathrm{(i)}\Rightarrow\mathrm{(ii)}$ follows from the
definition of a PCML matrix.

To prove $\mathrm{(ii)}\Rightarrow\mathrm{(iii)}$, suppose that
$\mathcal G_X$ contains a directed cycle $i_1\longrightarrow i_2\longrightarrow\cdots
\longrightarrow i_r\longrightarrow i_1.$ By the definition of $\mathcal G_X$, the submatrix consisting of these
columns, in the indicated order, is
\[
\begin{pmatrix}
x_{i_1}&x_{i_2}&\cdots&x_{i_r}\\
x_{i_2}&x_{i_3}&\cdots&x_{i_1}
\end{pmatrix}.
\]
The identity in the proof of Lemma~\ref{lem:app-cycle} yields $(x_{i_1}-x_{i_2})(x_{i_1}+\cdots+x_{i_r})\in I_2(X).$ Neither factor belongs to $I_2(X)$, because both are nonzero linear
forms. This contradicts $\mathrm{(ii)}$.

We next prove $\mathrm{(iii)}\Rightarrow\mathrm{(iv)}$. Since
$\mathcal G_X$ has no directed cycle, each of its connected components
is a directed path. Let $i_1\longrightarrow i_2\longrightarrow\cdots\longrightarrow i_r$ be such a component. Then $z_{i_k}=x_{i_{k+1}}$ for $1\le k<r$.
The terminal vertex $i_r$ has outdegree zero, so $z_{i_r}$ does not
occur in the first row. The columns belonging to this component, ordered
as $i_1,\ldots,i_r$, therefore form the block
\[
\begin{pmatrix}
x_{i_1}&x_{i_2}&\cdots&x_{i_r}\\
x_{i_2}&x_{i_3}&\cdots&z_{i_r}
\end{pmatrix}.
\]
Its $r+1$ variables are pairwise distinct. Moreover, blocks arising from
different components have disjoint sets of variables. Indeed, a variable
occurring in both blocks would have to occur in the first row of one
block and the second row of the other, and would therefore determine an
edge joining the two components. Relabeling the variables of each block
gives~\eqref{eq:app-scroll-blocks}. An isolated vertex gives a block with
$r_j=1$.

Finally, assume $\mathrm{(iv)}$. Write
$A=K[w_{j,k}:1\le j\le s,\ 0\le k\le r_j]$.
The matrix $H$ is the standard scroll matrix of type
$(r_1,\ldots,r_s)$. By
\cite[Appendix~A2H, p.~224]{EisenbudSyzygies},
$I_2(H)$ is prime. Moreover,
\cite[Corollary~A2.62]{EisenbudSyzygies}
shows that $I_2(H)$ has height $n-1$ and that
$A/I_2(H)$ is Cohen--Macaulay. The Eagon--Northcott complex of the
$2\times n$ matrix $H$ is therefore a minimal graded free resolution of
$A/I_2(H)$. Its module in homological degree $i\ge1$ is generated in
degree $i+1$, so $I_2(H)$ has a $2$-linear resolution; see
\cite[Appendix~A2H]{EisenbudSyzygies}. These properties are preserved
under the indicated permutation and relabeling. Hence $X$ is PCML.
\end{proof}

%=%=%=%=%=%=%=%=%=%=%=%=%=%=%=%=%=%=%=%=%=%=%=%=%=%=%=%=%=%=%=%=

%=%=%=%=%=%=%=%=%=%=%=%=%=%=%=%=%=%=%=%=%=%=%=%=%=%=%=%=%=%=%=%=
\begin{Remark}\label{rem:app-nearly-two-row}
{\em 
The class considered in this appendix allows variables in the second
row which do not occur in the first row. For a nearly determinantal
matrix $X_{2,n}$, every second-row entry occurs in the first row.
Since the entries in each row are pairwise distinct, every vertex of
$\mathcal G_{X_{2,n}}$ then has indegree and outdegree one. Thus the
graph is a disjoint union of directed cycles, all of length at least two
by~\textup{(ND2)}. Theorem~\ref{thm:app-two-row-classification} therefore
implies that $I(X_{2,n})$ is not prime.
}
\end{Remark}

%=%=%=%=%=%=%=%=%=%=%=%=%=%=%=%=%=%=%=%=%=%=%=%=%=%=%=%=%=%=%=%=%=%=%=%=

\bigskip

\begin{footnotesize}
{\bf Declaration of competing interest.}  The authors declare that they have no known competing financial interests or personal relationships that could have appeared to influence the work reported in this paper.

\medskip

{\bf Data availability.} No data was used for the research described in the article. 

\medskip

\end{footnotesize}

\end{document}